\documentclass[11pt]{amsart}

\usepackage{amsmath,amsthm, amscd, amssymb, amsfonts}
\usepackage[all]{xy}
\usepackage{verbatim}
\usepackage{color}
\usepackage{array}
\usepackage{multirow}
\usepackage{hhline}
\usepackage{enumitem}
\usepackage{float}
\allowdisplaybreaks
\usepackage{mathtools} 

\newtheorem{theorem}{Theorem}[section]
\newtheorem{lemma}[theorem]{Lemma}

\newtheorem{conjecture}[theorem]{Conjecture}
\newtheorem{proposition}[theorem]{Proposition}
\newtheorem*{claim}{Claim}

\theoremstyle{definition}

\newtheorem{example}[theorem]{Example}

\newtheorem*{notation}{Notation}

\theoremstyle{remark}
\newtheorem{remark}[theorem]{Remark}

\newcommand{\ydh}{{}^{H}_{H}\mathcal{YD}}

\def\pf{\begin{proof}}
\def\epf{\end{proof}}

\newcommand{\nc}{\newcommand}

\newcommand{\coloneqq}{:=}

\newcommand{\I}{{\mathbb I}}

\newcommand{\qti}{\widetilde{q}}

\nc{\ub}{\mathfrak{u}}
\nc{\g}{\mathfrak{g}}

\newcommand{\G}{{\mathbb G}}
\newcommand{\Gf}{{\mathbb G}_{\mathrm{f}}}

\newcommand{\Id}{\operatorname{Id}}

\newcommand\GL{\operatorname{GL}}
\newcommand\id{\operatorname{id}}

\newcommand\gdim{\operatorname{GKdim}}

\newcommand\gr{\operatorname{gr}}

\newcommand\co{\operatorname{co}}
\newcommand\ad{\operatorname{ad}}

\newcommand\rg{\operatorname{rg}}

\newcommand{\re}{\operatorname{re}}

\newcommand\ord{\operatorname{ord}}
\newcommand\Alg{\operatorname{Alg}}

\def\k{\Bbbk}

\def\ot{\otimes}

\def\C{\mathbb{C}}

\def\N{\mathbb{N}}
\def\Z{\mathbb{Z}}

\def\Re{\mathbb{R}}
\def\R{\mathcal{R}}
\def\Q{\mathbb{Q}}

\def\H{\mathcal{H}}

\def\mR{\mathcal{R}}

\def\mX{\mathcal{X}}

\def\eps{\varepsilon}

\newcommand\hgt{\operatorname{ht}}

\newcommand\Prim{\mathrm{P}}

\def\lg{\langle}
\def\rg{\rangle}

\newcommand{\bq}{{\bf q}}
\newcommand{\ba}{{\bf c}}
\newcommand{\bm}{{\bf m}}

\renewcommand{\tt}[1]{\texttt{#1}}
\newcommand{\qt}[1]{``#1''}
\newcommand{\one}{\tt{Crit1}}

\newcommand{\Ft}{\mathtt{F}}
\newcommand{\Gt}[1]{\mathtt{G}(#1)}
\newcommand{\Gtp}[1]{\mathtt{G}'(#1)}
\newcommand{\Gtpp}[1]{\mathtt{G}''(#1)}
\newcommand{\Ados}[1]{\mathtt{A}_2(#1)}
\newcommand{\Bdos}[1]{\mathtt{B}_2(#1)}
\newcommand{\Bdost}[1]{\mathtt{B}_2^t(#1)}
\newcommand{\Gdos}[1]{\mathtt{G}_2(#1)}
\newcommand{\Gdost}[1]{\mathtt{G}_2^t(#1)}
\newcommand{\sAuno}[1]{\mathtt{sA}_1(#1)}
\newcommand{\sAdos}[1]{\mathtt{sA}_2(#1)}
\newcommand{\sAdost}[1]{\mathtt{sA}_2^t(#1)}
\newcommand{\sB}[1]{\mathtt{sB}(#1)}
\newcommand{\sBt}[1]{\mathtt{sB}^t(#1)}
\newcommand{\Bst}[1]{\mathtt{Bst}(#1)}
\newcommand{\Bstt}[1]{\mathtt{Bst}^t(#1)}

\def\B{\mathfrak{B}}
\def\K{\mathcal{K}}

\def\gap{\texttt{GAP}}

\newcommand{\trianglerowone}{
\xymatrix@R-12pt{ \overset{1}{\circ} \ar  @{-}[rr]  &  &\overset{2}{\circ}\\
&    \underset{3}{\circ} \ar  @{-}[ur]\ar  @{=>}[ul] & 
}}

\newcommand{\Dtriangle}[6]{
\xymatrix@R-12pt{ \overset{#1}{\circ} \ar  @{-}[rr]^{#2}  &  &\overset{#3}{\circ}\\
 &    \underset{#5}{\circ} \ar  @{-}[ur]_{#4}\ar  @{-}[ul]^{#6} & 
 }}
\newcommand{\Dchaintwo}[3]{\xymatrix{\overset{#1}{\underset{\ }{\circ}}\ar
@{-}[r]^{#2}
& \overset{#3}{\underset{\ }{\circ}}}}

\newcommand{\lin}[3]{\xymatrix{\overset{#1}{\underset{\ }{\circ}}\ar
		@{-}[r]^{#2}
		& \overset{#3}{\underset{\ }{\circ}}}}

\newcommand{\liin}[5]{\xymatrix@C-4pt{
		\overset{#1}
		{\underset{\ }{\circ}}\ar@{-}[r]^{#2}
		& \overset{#3}{\underset{\ }{\circ}}\ar@{-}[r]^{#4}
		& \overset{#5}{\underset{\ }{\circ}}
	}
}

\newcommand{\Dchainthree}[5]{\xymatrix@C-4pt{
\overset{#1}
{\underset{\ }{\circ}}\ar@{-}[r]^{#2}
& \overset{#3}{\underset{\ }{\circ}}\ar@{-}[r]^{#4}
& \overset{#5}{\underset{\ }{\circ}}
}
}

\newcommand{\tri}[6]{
	\xymatrix@C-4pt@R-18pt{ & \overset{#5}{\underset{\ }{\circ}} \ar @{-}[ld]_{#6} \ar @{-}[rd]^{#4} & \\
		\overset{#1}{\underset{\ }{\circ}} \ar @{-}[rr]^{#2} & & \overset{#3}{\underset{\ }{\circ}}} }

\begin{document}



\title[Rank 3 and Cartan type]
{On finite GK-dimensional Nichols algebras of diagonal type: rank 3 and Cartan type}

\author[Angiono, Garc\'ia Iglesias]{Iv\'an Angiono, Agust\'in Garc\'ia Iglesias}

\address{FaMAF-CIEM (CONICET), Universidad Nacional de C\'ordoba,
Medina A\-llen\-de s/n, Ciudad Universitaria, 5000 C\' ordoba, Rep\'ublica Argentina.}

\email{(angiono|aigarcia)@famaf.unc.edu.ar}

\thanks{The work was partially supported by CONICET,
FONCyT-ANPCyT, Secyt (UNC), the MathAmSud project GR2HOPF}

\begin{abstract}
This paper contributes to the proof of the conjecture posed  in \cite{AAH-conj}, stating that a Nichols
algebra of diagonal type with finite Gelfand-Kirillov dimension has finite
(generalized) root system. We prove the conjecture assuming that the
rank is 3 or that the braiding is of Cartan type. 
\end{abstract}

\maketitle

\section{Introduction}\label{sec:intro}

This article picks up the task initiated in \cite{AAH-rank2}; that is to show that the Gelfand-Kirillov dimension of a Nichols algebra of diagonal type is finite if and only if its root system is finite. More precisely, as the converse is always true using the existence of a PBW basis indexed by the roots, the goal is to establish the validity of the following statement in \cite{AAH-conj}:

\begin{conjecture}
	\label{conjecture}
	Let $(V,c)$ be a braided vector space of diagonal type such that $\gdim \B(V)<\infty$. Then its generalized root system is finite.
\end{conjecture}
This conjecture is true when $\dim V=2$, by \cite[Theorem 1.2]{AAH-rank2}. In this article we show that it also holds when 
\begin{itemize}
	\item $V$ is of Cartan type, in Theorem \ref{thm:cartan}, or 
	\item $\dim V=3$, in Theorem \ref{thm:rank3}.
\end{itemize}

\medbreak

Nichols algebras $\B(V)$ with finite root systems were classified in \cite{H-rk2,H-rk3,H-full}, in terms of Dynkin diagrams representing the braiding matrix $\bq$ of $V$. This classification is presented in tables, according to the rank, that is, the dimension of $V$. Along our work, we refer to these tables as {\it Heckenberger's list}; or simply {\it the list}, for short. By abuse of notation, we will say that $V$, or the braiding matrix $\bq$, is in the list when the associated diagram is so. We write, equivalently, $\B_{\bq}$ to refer to $\B(V)$. 

The conjecture is thus equivalent to: 
\begin{center}
	{\it If  $\gdim \B_{\bq}<\infty$, then  $\bq$ is in the list.}
\end{center}

\medbreak

A positive answer to the conjecture has different implications towards the classification of pointed Hopf algebras with abelian coradical and finite $\gdim$. Indeed, the first step of Andruskiewitsch-Schneider \emph{Lifting Method} is to classify Nichols algebras over abelian groups with finite $\gdim$, and it gives the whole answer when the braiding is of diagonal type. Also, it is assumed to be true in \cite{AAH-conj} to obtain the classification of Nichols algebras of finite $\gdim$ whose braiding is a sum of \emph{points and blocks} (not of diagonal type), and a forthcoming paper on the general case of Nichols algebras over abelian groups. Finally, in \cite{ASa,ACS} the authors also assume the validity of the conjecture to classify pre and post-Nichols algebras of diagonal type with finite $\gdim$, which is part of the second step of the Lifting Method.

\medbreak

We remark that in Heckenberger's approach in \cite{H-rk3,H-full} to classify finite root systems, the cases of ranks three and four are key to attack the general case later on. We believe that this is also possible in our setting, and partially why we have chosen to focus on the rank-three case in this paper. Nevertheless, we also give a full answer for the braidings of Cartan type, of any rank. Again, this is also a crucial step in Heckenberger's work, and it is also endowed with a chronological importance, as these are the kind of braidings which have been consider first in many instances of the literature.

\medbreak

The structure of the paper is the following. In \S \ref{sec:pre} we recall some notions about Nichols algebras and root systems, focusing mainly on diagonal type: for the purposes of this article, we include a generalization of \cite[Proposition 3.1 \& Lemma 3.7]{AAH-rank2} to standard braidings. Next we deal with subquotients of Nichols algebras, which are obtained by means of the $\N_0^{\theta}$-graduation of Nichols algebras of diagonal type, see \S \ref{sec:subq}. Thanks to Proposition \ref{pro:subquotient} therein we will construct, for each $\omega\in\Z^{\theta}$, a new Nichols algebra generated by \emph{roots orthogonal to $\omega$}, a key step towards to check the validity of Conjecture \ref{conjecture} recursively on the rank; indeed we discuss the relationship between the root systems of our initial algebra and the new one obtained as subquotient, and apply this recursive machinery to give three criteria in \S \ref{sec:crit} to decide that some Nichols algebras have infinite $\gdim$ in the subsequent sections. Next, we deal in \S \ref{sec:cartan} with the case of Nichols algebras of Cartan type, and prove one of our main results in Theorem \ref{thm:cartan}, namely that a Nichols algebra of Cartan type has finite $\gdim$ if and only if the associated Cartan matrix is finite. The proof runs as follows: first we discuss how to reduce the problem to compactly hyperbolic matrices and then apply the criteria stated in \S \ref{sec:crit} to these matrices. Finally  in \S \ref{sec:rank3} we prove the validity of the conjecture for rank three in Theorem \ref{thm:rank3}. Again, we make use of the criteria developed in \S \ref{sec:crit}. Given a braiding matrix $\bq$ such that $\gdim \B_{\bq}<\infty$, one knows that each rank two submatrix has finite root system by \cite{AAH-rank2}; hence we take care of all possible ways to glue the rank two diagrams in the list and reduce to check just a finite  set in \S \ref{subsec:eval-param} using these criteria. This finite set of remaining diagrams is quite large, and we conclude the proof in \S \ref{subsec:gap} using computer program \gap\, where we implement the criteria to show that it contains matrices already on the list or such that $\gdim\B_{\bq}=\infty$, by our results on standard braidings from \S \ref{sec:pre}.

\section{Preliminaries}\label{sec:pre}
For $m\le n \in\N$, we denote by $\I_{m,n}=\{k\in\N: m\le k\le n\}$, $\I_n=\I_{1,n}$.

We work over an algebraically closed field of characteristic zero $\k$. 
Let $\G_N$ be the group of $N$th roots of 1 in $\k$ and let $\G'_N\subset \G_N$ be the subset of primitive roots; we also set $\G_{\infty}=\bigcup_{N\geq 2}\G_N$.

We recall some notation from quantum calculus. Consider the following expressions in the polynomial ring $\Z[t]$, for $n\geq k \geq 0$:
\begin{align*}
(n)_t &\coloneqq\sum_{j=0}^{n-1}t^j, & (n)_t!&\coloneqq\prod_{j=1}^{n}(j)_t, &
\binom{n}{k}_{\hspace*{-.1cm}t} \coloneqq \frac{(n)_t!}{(k)_t!(n-k)_t!}.
\end{align*}
If $q\in \k$; we denote by $(n)_q$, $(n)_q!$ and $\binom{n}{k}_{\hspace*{-.1cm}q}$ the corresponding evaluations.

\smallbreak

Let $A$ be a finitely generated $\k$-algebra and choose a finite-dimensional generating subspace $V\subseteq A$. Let $A_n=\sum_{j=0}^n V^j$. Then
\[
\gdim A =\limsup\limits_{n\to \infty} \log_n(\dim A_n).
\]
As expected, $\gdim A$ does not depend on the choice of $V$. We refer the reader to \cite{KL} for more information on this topic.

\smallbreak

We use Sweedler notation for the coproducts of Hopf algebras, or Hopf algebras in braided tensor categories. Given a Hopf algebra $H$, we denote by 
$\Prim(H)=\{x\in H: \Delta(x)=x\ot 1+1\ot x\}$ the space of primitive elements and $G(H)=\{x\in H-0: \Delta(x)=x\ot x\}$ the set of group-like elements.

\subsection{Nichols algebras of diagonal type}

Recall that a braided vector space $(V,c)$ is of diagonal type if there is a basis $\{x_1,\dots,x_\theta\}$ of $V$ for which the braiding is determined by a scalar matrix $\bq\in(\Bbbk^{\times})^{\theta\times\theta}$ as
\[
c(x_i\ot x_j)=q_{ij}x_j\ot x_i, \qquad i,j\in\I_\theta.
\]  
We refer to this $\bq$ as the braiding matrix of $(V,c)$ and we shall denote by $\B_{\bq}$ the corresponding Nichols algebra, which is a $\N_0$-graded connected Hopf algebra in a suitable category of Yetter-Drinfeld modules $\ydh$; namely, for any $H$ with a principal YD-realization $(\chi_i,g_i)_{i\in\I_\theta}\in \Alg(H,\k^\times)\times G(H)$. 

We refer to \cite{A-leyva,AA-diag-survey} for more details on Nichols algebras and realizations. 

\begin{notation}
	For any $x\in\B_{\bq}$ with $\deg x=n>1$ we have
	\[
	\Delta(x)=x\ot 1 + 1\ot x+\sum_{k=1}^{n-1} x_k\ot x^k, \qquad x_k,x^k\in\B_{\bq},
	\] 
	where $\deg x_k=k=n-\deg x^k$. We will use the following notation
	\[
	\underline{\Delta}(x):= \Delta(x)-x\ot 1 -1\ot x = \sum_{k=1}^{n-1} x_k\ot x^k,
	\]
	for the restricted comultiplication in $\B_{\bq}$.
\end{notation}

\subsection{Root systems}
The Nichols algebra $\B_{\bq}$ also inherits a $\Z^{\theta}$-grading, in such a way that there is a subset $L\subset \B_{\bq}$ of $\Z^{\theta}$-homogeneous elements together with a height function $\hgt\colon L\to \N\cup\{\infty\}$ so that
\begin{equation}\label{eqn:pbw}
\{\ell_1^{m_1}\dots \ell_k^{m_k} : k\in\N_0, \ell_1>\dots >\ell_k\in L, 0<m_i<\hgt(\ell_i), i\in\I_k \}
\end{equation}
is a linear basis of $\B_{\bq}$ \cite{K1}. The set of positive roots of $\bq$ is
\[
\varDelta^{\bq}_+=\{\deg \ell : \ell\in L  \}.
\]
By \cite[Lemma 4.7]{HS}, $\varDelta_+\coloneqq \varDelta_+^{\bq}$ does not depend on the choice of $L$. 

Moreover, let $K\subseteq \B_{\bq}$ be a right coideal subalgebra. Then $K$ admits a PBW basis as in \eqref{eqn:pbw} that can be extended to a basis of $\B_{\bq}$, see \cite{K2}.

As well, $\varDelta_+$ allows us to describe the (multivariate) Hilbert series associated to this graded algebra. Namely, for each $\alpha=n_1\alpha_1+\dots+n_\theta\alpha_{\theta}\in\varDelta_+$ we set $t^\alpha\coloneqq t_1^{n_1}\dots t_{\theta}^{n_{\theta}}\in \Z[[t_1,\dots,t_\theta]]$. Then the Hilbert series of  $\B_{\bq}$ can be factorized as:
\[
H_{\bq}(t)=\prod_{\alpha\in\varDelta_+}\varphi_{\hgt(\alpha)}(t^\alpha),
\]
where $\varphi_{\hgt(\alpha)}$ is the $\hgt(\alpha)$th cyclotomic polynomial.

\subsection{Weyl groupoid}\label{subsec:Weylgpd}

Next we recall the definition and properties of the Weyl groupoid of a Nichols algebra of diagonal type following the notation in \cite{AA-diag-survey}. We also compute roots of some examples needed in the sequel.

\medbreak
Fix $i\in \I$. We say that \emph{we can reflect $\bq$ at} $i$ if for all $j\neq i$ there exists $n\in\N_0$ such that $(n+1)_{q_{ii}}(1-q_{ii}^n q_{ij}q_{ji})=0$. If so, 
we set as in \cite{H-inv} the generalized Cartan matrix $C^{\bq}=(c_{ij}^{\bq})$, where $c_{ii}^{\bq} =2$ and
\begin{align}\label{eq:defcij}
c_{ij}^{\bq}&:=-\min\{n\in\N_0:(n+1)_{q_{ii}}(1-q_{ii}^nq_{ij}q_{ji})=0\},& & &j\neq i.
\end{align}
We also set $s_i^{\bq} \in \mathrm{GL}(\Z^{\theta})$ given by 
\begin{align}\label{eq:def-si}
s_i^{\bq} (\alpha_j) &= \alpha_j - c_{ij}^{\bq}\alpha_i, &  j & \in \I.
\end{align}

We use $s_i$ to define a new matrix $\rho_i \bq  = (t_{jk})_{j,k\in\I}$, called the reflection at the vertex $i$ of $\bq$, where
\begin{align}\label{eq:def-rho-ij}
t_{jk}&:= q_{s_i^{\bq}(\alpha_j), s_i^{\bq}(\alpha_k)} 
=  q_{jk}q_{ik}^{-c_{ij}^{\bq}}q_{ji}^{-c_{ik}^{\bq}}q_{ii}^{c_{ij}^{\bq} c_{ik}^{\bq}}, & j,k&\in\I.
\end{align}
Let $\rho_i V$ be the braided vector space of diagonal type with matrix $\rho_i \bq $. Notice that
$c_{ij}^{\rho_i \bq }=c_{ij}^{\bq}$ for all $j\in\Id$, so
\begin{align*}
s_{i}^{\rho_i \bq }&=s_{i}^{\bq}, & \rho_i (\rho_i \bq)  &=\bq.
\end{align*}
By \cite{H-inv, AA},
\begin{align}\label{eq:si-Delta}
\varDelta^{\rho_i \bq }_+ &= s_i^{\bq}(\varDelta^{\bq}_+-\{\alpha_i\})\cup \{\alpha_i\}, 
\\\label{eq:GKdim-reflection}
\gdim \B_{\rho_i \bq } &= \gdim \B_{\bq}.
\end{align}

We say that $\bq'$ is \emph{Weyl equivalent} to $\bq$ if there exist $i_j\in\I$ such that $\bq'=\rho_{i_k} \dots \rho_{i_1} \bq$.

We say that $\bq$ \emph{admits all reflections} if we can reflect $\bq'$ at every $i\in \I$ for all $\bq'$ Weyl equivalent to $\bq$.

If $\bq$ admits all reflections, then we denote by $\mX_{\bq}$ the collection of all braided vector spaces of diagonal type  obtained from $\bq$ by a finite number of successive reflections at various vertices (the equivalence class under the relation above). This happens for example when $\gdim \B_{\bq}<\infty$. 

The basic datum of $\bq$ is the pair $(\mX,\mR)$, where $\mR:\I\to\mathbb{S}_{\mX}$, $i\mapsto \rho_i$. The graph of the basic datum has $\mX$ as set of points, and an arrow between $\bq'$ and $\rho_i \bq'$ labelled with $i$ for each $\bq'\in\mX$ and $i\in\I$.
The \emph{Weyl groupoid} of $\bq$ is the subgroupoid of $\mX\times\GL(\Z^{\theta})\times \mX$ generated by
\begin{align*}
\sigma_i^{\bq'} & := (\bq', s_i^{\bq'},\rho_i\bq'), & & \bq'\in\mX, i\in\I.
\end{align*}
We use the following notation: For $i_j\in\I$,
\begin{align*}
\sigma_{i_1}^{\bq'}\sigma_{i_2}\cdots \sigma_{i_k}:=
\sigma_{i_1}^{\bq'}\sigma_{i_2}^{\rho_{i_1}\bq'} \cdots \sigma_{i_2}^{\rho_{i_{k-1}}\cdots \rho_{i_1}\bq'}.
\end{align*}
i.e., the implicit superscripts are the only possible allowing compositions. Similarly,
\begin{align*}
s_{i_1}^{\bq'}s_{i_2}\cdots s_{i_k}:=s_{i_1}^{\bq'}s_{i_2}^{\rho_{i_1}\bq'} \cdots s_{i_2}^{\rho_{i_{k-1}}\cdots \rho_{i_1}\bq'}.
\end{align*}

The collection $(\varDelta^{\bq'}_+)_{\bq' \in \mX}$ is the \emph{generalized root system} of $\bq$.

A root $\beta\in\varDelta^{\bq}$ is \emph{real} if there exist $i_j, i\in\I$ such that
\begin{align*}
\beta &= s_{i_1}^{\bq}s_{i_2}\cdots s_{i_k}(\alpha_i).
\end{align*}
We denote by $\varDelta^{\bq,\re}$ the set of all real roots, and $\varDelta^{\bq,\re}_+=\varDelta^{\bq,\re} \cap \N_0^{\I}$.

\begin{remark}\label{rem:enough}
	The roots and the reflections really depend on the twist equivalence class of $\bq$, that is, on the Dynkin diagram of $\bq$. Hence we will work with roots and reflections of equivalence classes of matrices (Dynkin diagrams) more than with matrices themselves.
	
	In particular, we can choose our matrix $\bq$ so that $q_{ij}=1$ for every $i>j$.
\end{remark}

\begin{notation}
	The vertices of a Dynkin diagram will be numbered from left to right and from bottom to top. For each $\tau\in\mathbb S_{\I}$, $\tau\bq:=(q_{\tau(i)\tau(j)})_{i,j\in\I}$. For $i\ne j\in\I$, $\tau_{ij}$ is the transposition interchanging $i$ and $j$.
\end{notation}

\begin{example}\label{ex:basic-datum-6diagrams}
	Let $\zeta\in\G_3'$, $\bq^{(1)} = \xymatrix@C-4pt{\overset{-1}{\underset{\ }{\circ}} \ar @{-}[r]^{-\zeta} & \overset{\zeta}{\underset{\ }{\circ}} \ar @{-}[r]^{-\zeta} & \overset{-\zeta^2}{\underset{\ }{\circ}} }$.
	
	The basic datum of $\bq^{(1)}$ is:
	\begin{center}
		\begin{tabular}{c c c c c}
			$\underset{\bq^{(1)}}{\bullet}$
			& \hspace{-5pt}\raisebox{3pt}{$\overset{1}{\rule{30pt}{0.5pt}}$}\hspace{-5pt}
			& $\underset{\bq^{(2)}}{\bullet}$
			& \hspace{-5pt}\raisebox{3pt}{$\overset{2}{\rule{30pt}{0.5pt}}$}\hspace{-5pt}
			$\underset{\bq^{(3)}}{\bullet}$
			& \hspace{-5pt}\raisebox{3pt}{$\overset{1}{\rule{30pt}{0.5pt}}$}\hspace{-5pt}
			$\underset{\tau_{12}(\bq^{(4)})}{\bullet}$
			\\
			& & & \hspace{18pt}{\scriptsize 3} \vline & 
			\\
			& && \hspace{18pt} $\underset{\tau_{23}(\bq^{(5)})}{\bullet}$
			& \hspace{-5pt}\raisebox{3pt}{$\overset{2}{\rule{30pt}{0.5pt}}$}\hspace{-5pt}
			$\underset{\tau_{23}(\bq^{(6)})}{\bullet}$
		\end{tabular}
	\end{center}
	where
	\begin{align*}
	\bq^{(2)} &= \xymatrix@C-4pt{\overset{-1}{\underset{\ }{\circ}} \ar @{-}[r]^{-\zeta^2} &
		\overset{\zeta^2}{\underset{\ }{\circ}} \ar @{-}[r]^{-\zeta} & \overset{-\zeta^2}{\underset{\ }{\circ}} },
	\\
	\bq^{(3)} &= \xymatrix@C-4pt@R-18pt{ & \overset{-1}{\underset{\ }{\circ}} \ar @{-}[ld]_{\zeta} \ar @{-}[rd]^{-1} & \\
		\overset{-1}{\underset{\ }{\circ}} \ar @{-}[rr]^{-\zeta^2} & & \overset{\zeta^2}{\underset{\ }{\circ}}},
	&
	\bq^{(4)} &= \xymatrix@C-4pt{\overset{\zeta}{\underset{\ }{\circ}} \ar @{-}[r]^{-\zeta} &
		\overset{-1}{\underset{\ }{\circ}} \ar @{-}[r]^{\zeta^2} & \overset{\zeta}{\underset{\ }{\circ}} },
	\\
	\bq^{(5)} &= \xymatrix@C-4pt{\overset{\zeta}{\underset{\ }{\circ}} \ar @{-}[r]^{\zeta^2} &
		\overset{-1}{\underset{\ }{\circ}} \ar @{-}[r]^{-1} & \overset{\zeta^2}{\underset{\ }{\circ}} },
	&
	\bq^{(6)} &= \xymatrix@C-4pt{\overset{\zeta}{\underset{\ }{\circ}} \ar @{-}[r]^{\zeta^2} &
		\overset{-\zeta^2}{\underset{\ }{\circ}} \ar @{-}[r]^{-\zeta} & \overset{\zeta^2}{\underset{\ }{\circ}} }.
	\end{align*}
\end{example}

\subsection{Standard and Cartan braidings}\label{subsec:cartan-standard}
Assume that $\bq$ admits all reflections.
We say that $\bq$ is \emph{standard} if $C^{\bq'}=C^{\bq}$ for all $\bq'$ Weyl equivalent to $\bq$ \cite{AA}.
Also, $\bq$ is of \emph{Cartan} type if $q_{ii}^{c_{ij}^{\bq}}=\qti_{ij}$ for all $i\ne j\in\I$.
One can check that every braiding of Cartan type is standard; more generally, every $\bq$ such that $\rho_i \bq$ has the same Dynkin diagram as $\bq$ for all $i\in\I$ is standard.

\begin{remark}\label{rem:standard-roots}
	Let $\bq$ be a standard braiding with Cartan matrix $C=C^{\bq}$. Let $W$ be the Weyl group of $C$. By \eqref{eq:si-Delta}, $W\gamma \subseteq \varDelta^{\bq}$ for all $\gamma\in\varDelta^{\bq}$. 
	
	In particular, if $\varDelta^{\re}$ denotes the set of real roots of $C$, then
	$$ \varDelta^{\re} \subseteq \varDelta^{\bq, \re} \subseteq \varDelta^{\bq} .$$
\end{remark}

Next we extend two results from \cite{AAH-rank2} stated there for Cartan type to standard type: although the proofs are essentially the same, we repeat them here for completeness.

\begin{proposition}\label{prop:standard-affine}
	Let $\bq$ be a braiding of standard type with affine Cartan matrix $C=C^{\bq}$. Then $\gdim \B_{\bq}=\infty$.
\end{proposition}
\pf
As in the proof of \cite[Proposition 3.1]{AAH-rank2}, let $\varDelta^{\re}$ be the set of real roots of $C$: By \cite[Proposition 6.3 d)]{K} there exists a positive root $\delta$ such that 
$\varDelta^{\re}+\delta=\varDelta^{\re}$. Fix a homogeneous restricted PBW basis of $\B_{\bq}$. 
If $m$ is the height of $\delta$ and $\alpha$ is a simple root, Remark \ref{rem:standard-roots} says that there exists a PBW-generator of degree $k\delta+\alpha$ for all $k\ge 0$. Hence $\gdim \B_{\bq}=\infty$ by \cite[Lemma 2.3.4]{AAH-conj}.
\epf

\begin{lemma}\label{lem:standard-indefinite}
	Let $\bq$ be a braiding of standard type with $C=C^{\bq}$ indefinite. If there exists $\gamma\in\varDelta_+^{\bq}$ such that $q_{\gamma\gamma}=1$, then $\gdim \B_{\bq}=\infty$.
\end{lemma}
\pf
We follow the same ideas as for \cite[Lemma 3.7]{AAH-rank2}. By Remark \ref{rem:standard-roots}, $W\gamma\cap \N_0^{\I}\subseteq \varDelta^{\bq}_+$, and by \cite[Lemma 3.6]{AAH-rank2}, $W\gamma \cap \N_0^{\I}$ is infinite. Hence any homogeneous restricted PBW basis of $\B_{\bq}$ has infinite PBW generators of infinite height (one for each element in $W\gamma \cap \N_0^{\I}$). By \cite[Lemma 2.3]{AAH-rank2}, $\gdim\B_{\bq}=\infty$.
\epf

\section{Subquotients and root systems}\label{sec:subq}

We fix a braided vector space $(V,c)$ of diagonal type with braiding matrix $\bq$, $\B_{\bq}=\B(V)$ the corresponding Nichols algebra, and a Hopf algebra $H$ with bijective antipode such that there is a principal realization $V\in\ydh$.

\subsection{Subquotients of Nichols algebras}
Let $\B$ be a pre-Nichols algebra over $V$, that is an intermediate $\N_0$-graded quotient $T(V)\twoheadrightarrow \B\twoheadrightarrow \B_{\bq}$ in $\ydh$.
We shall further assume that $\B$ is a $\Z^{\theta}$-graded Hopf algebra in $\ydh$,
$\B=\bigoplus_{\alpha\in\Z^{\theta}}\B^\alpha$. We denote by $\pi_\alpha\colon \B\to \B^\alpha$ the linear projection onto the $\alpha$-component. Let $\alpha_1,\dots,\alpha_\theta$ be the canonical basis of $\Z^{\theta}$ and let $(\cdot|\cdot)$ be the inner product in $\Re^{\theta}$ defined by $(\alpha_i|\alpha_j)=\delta_{i,j}$. 

In this setting, for a fixed $\omega\in \Z^{\theta}$ we set
\begin{align*}
\B_{\geq 0}&\coloneqq \bigoplus_{\alpha:(\alpha|\omega)\geq 0} \B^\alpha, &  \B_{> 0}&\coloneqq \bigoplus_{\alpha:(\alpha|\omega)> 0} \B^\alpha.
\end{align*}
Clearly $1\in  \B_{\geq 0}$ and $\B_{> 0}\subset \B_{\geq 0}$. We also set $\B_{0}\coloneqq \B_{\geq 0}\setminus \B_{> 0}$ and $\B_{<0}$ as expected.
As we well, we define 
\begin{align*}
\K_{\geq 0}&\coloneqq \{x\in \B:\Delta(x)\in \B_{\geq 0}\ot \B\}, &  \K_{> 0}&\coloneqq \K_{\geq 0}\cap \B_{>0}
\end{align*}
and $\K_{0}\coloneqq \K_{\geq 0}\cap \B_{0}$, $\K_{<0}$, accordingly. We will make use of the induced linear decompositions $\B=\B_{<0}\oplus \B_{0}\oplus \B_{>0}$ and $\K_{<0}\oplus \K_{0}\oplus \K_{>0}\subset \B$.

\begin{remark}
	Lemma \ref{lem:K} and Proposition \ref{pro:subquotient} next extend \cite[Proposition 3.9]{AAH-rank2}, where $\theta=2$ and $\omega=\alpha_1-r\,\alpha_2$, $0\leq r\in \Q$. The proofs use the same ideas as in loc.cit.~but nevertheless we choose to include a detailed proof for the sake of completeness.
\end{remark}

\begin{lemma}\label{lem:K}
	For $\B$ and $\omega\in\Z^{\theta}$ as above, we have that $\B_{\geq 0}\subset \B$ is a subalgebra in $\ydh$. Besides,
	\begin{enumerate}
		\item $\K_{\geq 0}\subseteq \B_{\geq 0}$.
		\item $\K_{\geq 0}\subset \B$ is a subalgebra in $\ydh$.
		\item $\Delta(\K_{\geq 0})\subseteq \K_{\geq 0}\ot \K_{\geq 0}+\K_{> 0}\ot \B$.
		\item $\K_{>0}$ is an ideal of $\K_{\geq 0}$ and a coideal of $\B$ in $\ydh$.
	\end{enumerate}
\end{lemma}
\pf
If $x\in \B^\alpha\cap \B_{\geq 0}$ and $y\in \B^\beta\cap \B_{\geq 0}$, then $xy\in \B_{\alpha+\beta}$ and $(\alpha+\beta|\omega)=(\alpha|\omega)+(\beta|\omega)\geq 0$. Hence $\B_{\geq 0}$ is a subalgebra of $\B$. Now, if $x\in \K_{\geq 0}$, then $x=x_{(1)}\eps(x_{(2)})\in \B_{\geq 0}$ by definition, so (1) follows. If $x,y \in \K_{\geq 0}$, then $xy \in \K_{\geq 0}$ since the comultiplication is an algebra map.

Now fix $x\in \B^\alpha\cap \K_{\geq 0}$ and let us write $\Delta(x)=\sum y_i\ot z_i\in \B_{\geq 0}\ot \B$, with $y_i$, resp.~$z_i$, homogeneous of degree $\beta_i$, resp.~$\gamma_i=\alpha-\beta_i$. On the one hand, 
\[
\sum \Delta(y_i)\ot z_i=\sum y_i\ot \Delta(z_i)\in \B_{\geq 0}\ot \B\ot \B,
\]
which shows that $\Delta(\K_{\geq 0})\subseteq \K_{\geq 0}\ot \B$, i.e.~it is a right coideal. Now, if $(\beta_i|\omega)>0$ for some $i$, then $y_i\ot z_i\in \K_{>0}\ot \B$. If $(\beta_i|\omega)=0$, then $(\alpha-\beta_i,\omega)\geq0$ and thus $y_i\ot z_i\in \K_{\geq 0}\ot \B_{\geq 0}$. That is, 
\begin{align}\label{eqn:DeltaK1}
\Delta(\K_{\geq 0})\subseteq \K_{\geq 0}\ot \B_{\geq 0}+\K_{>0}\ot \B.
\end{align}
We have actually shown that $\Delta(\K_{\geq 0})\subseteq \K_{0}\ot \B_{\geq 0}+\K_{>0}\ot \B$.
Observe that this also yields $\Delta(\K_{> 0})\subseteq \K_{\geq 0}\ot \B_{\geq 0}+\K_{>0}\ot \B$, which can actually be more accurately written as
\begin{align}\label{eqn:DeltaK2}
\Delta(\K_{> 0})\subseteq \K_{0}\ot \B_{> 0}+\K_{>0}\ot \B,
\end{align}
with the same argument as above.

Now, for $x$ as above, we use \eqref{eqn:DeltaK1} to write $\Delta(x)=\sum\limits_{i\in I} a_i\ot b_i+\sum\limits_{j\in J} c_j\ot d_j$, where $a_i\in \K_{0}$, $b_i\in \B_{\geq 0}$ and $c_j\in \K_{>0}$. To prove (3), we need to show that $b_i\in \K_{\geq 0}\subset \B_{\geq 0}$. Assume that, for some $i\in I$, $b_i\notin \K_{\geq 0}$, which is to say that $\Delta(b_i)\in \B_{<0}\ot \B$. In other words,
\begin{align*}
(\id\ot\Delta)\Delta(x)=\sum_i a_i\ot \Delta(b_i)+\sum_j c_j\ot\Delta(d_j)
\end{align*}
linearly projects non-trivially onto $\K_{0}\ot \B_{<0}\ot \B$. On the other hand, we can use \eqref{eqn:DeltaK1} and \eqref{eqn:DeltaK2} to deduce that
\begin{multline*}
(\Delta\ot\id)\Delta(x)\in (\Delta\ot\id)\Delta(\K_{\geq 0})\subset \Delta(\K_{\geq 0})\ot \B_{\geq 0}+\Delta(\K_{>0})\ot \B\\
\subset \K_{\geq 0}\ot \B_{\geq 0}\ot \B_{\geq 0}+\K_{>0}\ot \B\ot \B_{\geq 0}\\+\K_{0}\ot \B_{> 0}\ot \B+\K_{>0}\ot \B\ot \B
\end{multline*}
which does not intersect the subspace $\K_{0}\ot \B_{<0}\ot \B$. Thus $b_i\in \K_{\geq 0}$ for all $i\in I$ and (3) follows.

This also shows that, in particular, $\Delta(\K_{>0})\subseteq \K_{\geq 0}\ot \K_{\geq 0}+\K_{> 0}\ot \B$. If $(\beta_i|\omega)=0$, then necessarily $(\alpha-\beta_i|\omega)=(\alpha|\omega)> 0$, so $\Delta(\K_{>0})\subseteq \K_{\geq 0}\ot \K_{>0}+\K_{> 0}\ot \B$ and thus $\K_{>0}$ is a coideal of $\B$. Finally, $\K_{>0}$ is an ideal of $\K_{\geq 0}$ as the multiplication is $\Z^{\theta}$-graded, which gives (4).
\epf

\begin{proposition}\label{pro:subquotient}
	Fix $\omega\in\Z^{\theta}$. Then the braided bialgebra structure of $\B$ induces a braided bialgebra structure on 
	\begin{equation*}
	\K_{\omega}\coloneqq \K_{\geq 0}/\K_{>0}.
	\end{equation*}
\end{proposition}
\pf
This is \cite[Proposition 3.8]{AAH-rank2}, using Lemma \ref{lem:K}.
\epf	

Notice that $\K_{\omega}\in\ydh$ and thus $\H_\omega =\K_{\omega}\#H$ is a Hopf algebra. 
The grading on $\K_\omega$ is not induced by the coradical filtration on $\H_\omega$ which in turn determines a graded Hopf algebra $\R_\omega=\oplus_{n\geq 0} \R_\omega^n\in\ydh$  so that $\gr\H_\omega=\R_\omega\# H$. We set
$V'=\R_\omega^1$; and denote by $\bar\bq$ its associated braiding matrix. In a snapshot, we have the following situation:
\begin{equation}\label{eqn:snapshot}
\xymatrix{
	K_{> 0}\ar@{^{(}->}[r] & K_{\geq 0}\ar@{^{(}->}[r]\ar@{->>}[d]       & \B_{\bq}\in \ydh  \\
	& \K_{\omega}\ar@{~>}^{\gr (-\#H)^{\co H}}[d]                 & E\coloneqq\k\lg x_\alpha, x_\beta\rg \ar@{_{(}->}[l]\\
	& \R & \B_{\bar\bq}\coloneqq \B(V')\ar@{_{(}->}[l]}
\end{equation}
Notice that $\gdim\B_{\bar\bq}=\infty\Rightarrow\gdim \B_{\bq}=\infty$. Based on this picture, we will develop in \S \ref{sec:crit} a series of criteria to show that a given matrix $\bar\bq$ has $\gdim \B_{\bq}=\infty$, by picking the suitable $\omega$. We remark that these ideas apply for braided vector spaces of any rank.

\subsection{An example}

We include an example of this strategy, for a diagram (and its reflections) that could not be tackled down with the criteria. 
This represents the unique case that escaped the criteria, and for which we had to apply an ad-hoc analysis, also based in Proposition \ref{pro:subquotient}.

\begin{lemma}\label{lem:grupoide6puntos}
	Let $\bq$ be any of the braiding matrices in Example \ref{ex:basic-datum-6diagrams}.
	Then $\gdim \B_{\bq}=\infty$.
\end{lemma}
\pf
By \eqref{eq:GKdim-reflection} it is enough to prove that $\gdim \B_{\bq}=\infty$ for one of the matrices in Example \ref{ex:basic-datum-6diagrams}. 
We choose $\bq=\bq^{(3)}$ and fix $\bq=
\left(\begin{smallmatrix}
-1&-\zeta^2&\zeta\\
1 & \zeta^2  &-1\\
1 &  1  &-1
\end{smallmatrix}\right)$, $\zeta\in\G_3'$. 
We shall consider $\omega=3\alpha_1-\alpha_2-\alpha_3$, and $\alpha=\alpha_1+2\alpha_2+\alpha_3$, $\beta=\alpha_1+\alpha_2+2\alpha_3$; so that $(\alpha,\omega)=(\beta,\omega)=0$.

We let $x_\alpha=[[x_{13},x_2]_c,x_2]_c$ and $x_\beta=[x_{13},x_{23}]_c$; by direct computation,
\begin{align*}
\underline{\Delta}&(x_\alpha)=   3x_{13}\ot x_2^2 + (1-\zeta)x_{122}\ot x_3  \\
&+ (\zeta-\zeta^2)([x_{12},x_3]_c - x_{123} -2\zeta^2x_2x_{13} )\ot x_2 \\
&+ (1-\zeta)x_1\ot x_{223}  +3\zeta x_1\ot x_{23}x_2 + 3(1-\zeta)x_1\ot x_3x_2^2\\
&+3\zeta^2x_{12}\ot x_3x_2 + (1-\zeta^2)x_{12}\ot x_{23}, 
\\
\underline{\Delta}&(x_\beta)= (-2\zeta x_2x_{13} -2\zeta^2 x_{23}x_1  - 2[x_{12},x_3] + 2\zeta^2x_{23}x_1 -\zeta^2x_{123})\ot x_3.
\end{align*}
It follows that $x_\alpha,x_\beta\neq 0$ and $x_\alpha,x_\beta\in\K_{\geq 0}$. Furthermore, (the classes of) $x_\alpha$ and $x_\beta$ are primitive elements in $\K_{\geq 0}/\K_{>0}$.

Now, the diagram associated to $\k\{x_\alpha,x_\beta\}$ is $\Dchaintwo{\zeta}{\zeta^2}{1}$, which is not of finite type. Hence $\gdim \B_{\bq}=\infty$ and the Lemma follows.
\epf

\subsection{More on $\K_{\omega}$}
Next we make some remarks on the structure of the bialgebras $\K_{\omega}$.
We aim to understand the relation between $\K_{\omega}$ and $\B_{\bar\bq}$; as well as between their root systems. We believe it could be of interest to determine under which conditions one can assure that $\K_\omega$ is generated by $\Prim(\K_\omega)$. Similarly, to pin necessary and sufficient conditions so that one or both of the inclusions in \eqref{eqn:inclusions} below hold.

We start with a description of $\Prim(\K_\omega)$.
\begin{remark}\label{rem:primitives}
	Let $x\in\K_{\geq 0}$ be an homogeneous element, of degree $\alpha\in\Delta$, so that  $\bar x\in \K_\omega$ is primitive. Clearly, $(\alpha|\omega)=0$. Now 
	$\underline\Delta(x)\in \sum_\beta \B_{\bq}^\beta\ot \B_{\bq}^{\alpha-\beta}$, with $(\beta|\omega)\geq 0$, by definition of $\K_{\geq0}$; hence $(\gamma|\omega)\leq 0$. 
	Since $\bar x $ is primitive, we see that we necessarily have $(\beta|\omega)>0$. This shows that
	\begin{align*}
	\Prim(\K_\omega)=\{x_\alpha : & (\alpha|\omega)=0\text{ and }\alpha\text{ simple } \\
	\text{ or }&(\beta|\omega)(\alpha-\beta|\omega)<0 \text{ for all }\beta : (\pi_\beta\ot\id)\underline\Delta(x_\alpha)\neq 0\}.
	\end{align*}
\end{remark}

For the analysis of the root system, we need some notation.
\begin{notation}
	We shall write $\Delta_\omega$ for the root system associated to $\K_{\omega}$; while $\Delta'$ will denote the root system corresponding to the Nichols algebra 
	$\B(V')$.
\end{notation}

We will show in the next two lemmas that there is a chain of inclusions
\begin{equation}\label{eqn:inclusions}
\Delta'\subseteq \Delta_{\omega}\subseteq \Delta\cap\{\omega\}^{\perp}.
\end{equation}

We begin with the first inclusion in \eqref{eqn:inclusions}.
\begin{lemma}
	$\Delta'\subseteq \Delta_{\omega}$.
\end{lemma}
\pf
To start with, we shall assume, without lost of generality, that $\bar\bq$ is symmetric, so $\B_{\bar\bq}^\ast\simeq \B_{{}^t\bar\bq}=\B_{\bar\bq}$ and thus the inclusion $\B_{\bar\bq}\hookrightarrow \R$ produces a projection $\R^\ast\twoheadrightarrow \B_{\bar\bq}$. As well recall that the Hilbert series $H_{\R}(t)$ and $H_{\R^\ast}(t)$ of  $\R$ and $\R^\ast$ coincide. Hence we have a factorization $H_{\R}(t)=H_{\bar\bq}(t)P(t)$, for some multivariate series $P(t)$.  This shows that
$\Delta'\subseteq \Delta_{\mR}=\Delta_{\omega}$.
\epf

As for the second, we have the following.
\begin{lemma}
	$\Delta_{\omega}\subseteq \Delta\cap\{\omega\}^{\perp}$.
\end{lemma}
\pf
By \cite{K2}, there is a set $L_{\geq 0}\subset K_{\geq 0}$ of $\Z^\theta$-homogeneous $PBW$ generators for  $\K_{\geq 0}$ which can be extended to a set of PBW generators $L\subset \B_{\bq}$. Thus $\Delta_{\omega}\subseteq \Delta$. 
Notice that $(\deg\ell|\omega)\geq 0$ for any $\ell\in L_{\geq 0}$, since $\K_{\geq 0}\subset \B_{\geq 0}$. 

Let us write $x\mapsto \bar{x}$ for the projection $K_{\geq 0}\twoheadrightarrow \K_{\omega}$ and set
\begin{align*}
L_0&=\{\ell\in L_{\geq 0}:(\deg\ell|\omega)=0\}, & S&=\{\bar{\ell}:\ell\in L_0\}.
\end{align*}
Notice that $L_0$ may not be the whole subset of $L$ of $\omega$-orthogonal degrees.

We claim that $S$ is a set of PBW generators for $\K_{\omega}$, proving the lemma.

\smallbreak

Let $\bar{x}\in\K_{\omega}$ be homogeneous of degree $\alpha$; necessarily $(\alpha|\omega)=0$. Now $x\in K_{\geq0}$ is a linear combination of products of elements in $L'$, that is there $(\lambda_{i_1,\dots,i_k})\in\k$ so that we can write $x=\sum \lambda_{i_1,\dots,i_k}\ell_{i_1}\dots \ell_{i_k}$, with $\sum \beta_{i_j}=\alpha$ and $(\beta_{i_j}|\omega)\geq 0$. Now, if any $i_j$, some $i,j$, is such that $(\beta_{i_j}|\omega)> 0$, then $\ell_{i_1}\dots \ell_{i_{k}}\in K_{>0}$ and $\overline{\ell_{i_1}\dots \ell_{i_{k}}}=0$; so we can assume that every $i_j$ satisfies $(\beta_{i_j}|\omega)= 0$ and hence $S$ spans this quotient algebra. 

On the other hand, a nonzero linear combination $\sum \lambda_{i_1,\dots,i_k}\ell_{i_1}\dots \ell_{i_k}$ of products of elements in $L_{0}$ cannot  land on $K_{>0}$, which shows that the set $S$ is linearly independent, and thus it is a basis of $\K_{\omega}$.
\epf

\section{Criteria}\label{sec:crit}

Let $\bq$ a braiding matrix with diagram $d$. If  $\gdim\B_{\bq}<\infty$, then each subdiagram of $d$ with two vertices belongs to Heckenberger's list by \cite{AAH-rank2}. Moreover, for any $\omega$ and any two primitive elements $\tt{x}$ and $\tt{y}$ in $\K_{\omega}$, the Dynkin diagram $d'=d(\tt{x},\tt{y})$ associated to the braided vector space $\{\tt{x},\tt{y}\}$ also belongs to the list. 

Let $d$ be any diagram with three vertices obtained by pasting two (a line) or three (a triangle) diagrams with two vertices. We develop three criteria, based on three different families of vectors $\omega$. We apply these criteria to all such diagrams $d$, and also apply it to some of their reflections. For each diagram $d$, if we find a vector $\omega$ in one of the families (and, if necessary, a reflection) for which $d'$ is not in the list, then automatically this yields $\gdim\B_{\bq}=\infty$ and we can discard $d$.

\subsection{Criterium 1}\label{crit1}
To consider the rank 2 diagrams associated to the couples 
$(\alpha_i,\ell\alpha_j+\alpha_k)$, $1\leq \ell\leq m_{jk}$, $i,j,k$ all different, for 
\[
\omega=\alpha_j-\ell\,\alpha_k.
\]

Assume $j<k$, set $\mu_t=\prod_{s=0}^{t-1}(1-q_{jj}^s\widetilde{q}_{jk})$, $t\geq 0$. 
Then
\begin{align}\label{eqn:del112}
\underline{\Delta}((\ad_c x_j)^\ell(x_k))&=\sum_{t=0}^{\ell-1} \binom{n}{t}_{\hspace*{-.2cm}q_{jj}}\frac{\mu_\ell}{\mu_t}\, x_j^{\ell-t}\ot (\ad_c x_j)^t(x_k).
\end{align}
Hence $\tt{x}=x_i$ and $\tt{y}=(\ad_c x_j)^\ell(x_k)$ become primitive in the corresponding subquotient $\K_{\omega}$, as $((\ell-t)\alpha_j|\omega)=\ell-t\geq 1$, $0\leq t\leq \ell-1$.

On the other hand, when $j>k$, we consider 
\[
\tt{y}={\ad'}_c^\ell(x_j)(x_k) := [[[[x_{kj},x_j],x_j],\dots ],x_j]
\]
and there are $(a_t)\in\k$, $t=0,\dots,k-1$, such that
\begin{align}\label{eqn:del221}
\underline{\Delta}({\ad'}_c^\ell(x_j)(x_k))&=\sum_{s=0}^{\ell-1} a_t {\ad'}_c^t(x_j)(x_k)\ot x_j^{\ell-t}.
\end{align}
Hence $\tt{x}=x_i$ and $\tt{y}=(\ad_c x_j)^\ell(x_k)$ are primitive in $\K_{\omega}$.
%

\subsection{Criterium 2}\label{crit2}

To consider the rank 2 diagrams associated to the couples $(\alpha_1+\alpha_2+\alpha_3,\ell\alpha_i+\alpha_j)$, $2\leq \ell\leq m_{ij}$, $1\leq i\neq j\leq 3$, so
\[
\omega=c_1\alpha_1+c_2\alpha_2+c_3\alpha_3, \qquad  c_j=-kc_i, \  c_3=-(c_1+c_2).  
\]
We analyze all possible $(c_1,c_2,c_3)$, depending on the values of $i,j$.

$\diamondsuit$ Assume that $\qti_{12}, \qti_{23}\ne 1$, $\qti_{13}=1$, so
\begin{align*}
\underline{\Delta}(x_{123})&=q_{12}(1-\widetilde{q}_{23})(1-\widetilde{q}_{13}) x_2x_1\ot x_3 + (1-\widetilde{q}_{12}\widetilde{q}_{13})x_1\ot x_{23} \\
&\qquad + (1-\widetilde{q}_{23})x_{12}\ot x_3 + q_{12}(1- \widetilde{q}_{23})x_2\ot x_{13}\\
&=(1-\widetilde{q}_{12})x_1\ot x_{23} + (1-\widetilde{q}_{23})x_{12}\ot x_3.
\end{align*}	
%
Then $x_{123}\in \K_{\geq 0}$ when $c_1\geq 0$ and $c_2\geq -c_1$; while the (the class of) $\tt{x}=x_{123}$ is primitive in $\K_{\omega}$ when $c_1 > 0$ and $c_2>-c_1$. On the other hand, $\tt{y}={\ad'}^\ell(x_2)(x_1)$ is also primitive when $c_1=-\ell\,c_2$ that is
\[
\omega=\ell\alpha_1-\alpha_2+(\ell-1)\alpha_3.\]
To deal with $\ell\alpha_1+\alpha_2$ we consider $\tt{x}=[x_3,[x_2,x_1]_c]_c$. Similar choices, {\it mutatis mutandis}, are done for $\ell\alpha_2+\alpha_3$ and $\ell\alpha_3+\alpha_2$.

$\diamondsuit$
Now, if $\qti_{12}, \qti_{23}, \qti_{13}\ne 1$. In this case, we turn to $\tt{x}=[x_{13},x_2]_c$. As
\begin{align*}
\underline{\Delta}&([x_{13},x_2]_c)=q_{32}(1-\widetilde{q}_{13})x_{12}\ot x_3 +(1-\widetilde{q}_{12}\widetilde{q}_{23})x_{13}\ot x_2\\
&
-\widetilde{q}_{12}q_{32}(1-\widetilde{q}_{13})x_{1}\ot x_{23}+(1-\widetilde{q}_{12}\widetilde{q}_{23})(1-\widetilde{q}_{13})x_{1}\ot x_3x_{2},
\end{align*}
$\tt{x}$ is primitive in $\K_{\omega}$ when $-c_1 <c_2<0$, as this gives $c_1>0$, $c_1+c_2>0$ and $c_1+c_3=-c_2>0$.
As well, $\tt{y}={\ad'}^\ell(x_2)(x_1)$, $\ell>1$ is primitive in $\K_{\omega}$ for
\[
\omega=\ell\alpha_1-\alpha_2-\alpha_3.
\]
Similar choices are made for different pairs $i\neq j\in\I_3$.

\subsection{Criterium 3}\label{crit3}
To consider the rank 2 diagrams associated to the couples $(\alpha_i+\alpha_k,\alpha_j+\alpha_k)$, $i,j,k\in\{1,2,3\}$ all different, so
\[
\omega=c_i\alpha_i+c_j\alpha_j+c_k\alpha_k, \qquad  c_k=-c_i=-c_j.  
\]
This is clear. This criterium can be extended further by considering couples of the form $(n\alpha_i+\alpha_k,\alpha_j+\alpha_k)$, with $n\leq m_{ik}$. In this case, \[
\omega=\alpha_i-n\alpha_k+n\alpha_j.
\]
Notice that for  $\tt{x}=\ad(x_i)^n(x_k)$, the right choice of $\tt{y}=x_{jk}$ or $\tt{y}=x_{kj}$ has to be made (independent of the order of $k,j$). For instance, one may choose $\tt{x}=x_{112}$ and $\tt{y}=x_{32}$, as $x_{23}\notin \K_{>0}$, for $\omega=\alpha_1-2\alpha_2+2\alpha_3$.

\section{Cartan type}\label{sec:cartan}

Let $(V,c)$ be a braided vector space of Cartan type, see \S \ref{subsec:cartan-standard}, with matrix $\bq=(q_{ij})_{i,j\in\I_\theta}$. Let $\ba=(c_{ij}^{\bq})_{i,j\in\I_\theta}$, cf.~\eqref{eq:defcij},  be the associated generalized Cartan matrix (GCM for short). Indecomposable GCM are of three types: finite, affine or indefinite; we call the corresponding connected braided vector spaces of Cartan type and the matrix $\bq$ finite, affine o indefinite, accordingly.

We recall that a braiding matrix $\bq=(q_{ij})$ is called
\begin{itemize}[leftmargin=*]
	\item {\it generic} when $q_{ii}\notin\G_{\infty}$ and either $\widetilde{q}_{ij}=1$ or $\widetilde{q}_{ij}\notin\G_{\infty}$, for every $i,j\in\I_\theta$;
	\item {\it of torsion class} when $q_{ii},\widetilde{q}_{ij}\in\G_{\infty}$, for every $i,j\in\I_\theta$;
	\item {\it semigeneric} when it is neither generic nor of torsion class.
\end{itemize}

The next theorem states an equivalence between finite type and finite GK-dimension for $\B_{\bq}$.

\begin{theorem}\label{thm:cartan}
	Let $\bq$ be a matrix of Cartan type such that $\gdim \B_{\bq}<\infty$. Then its root system 
	is finite.
\end{theorem}

Before the proof, let us recall that an indefinite Cartan matrix is compactly hyperbolic if every proper minor is of finite type. These matrices have been classified in \cite{C}; we shall use such classification in our proof.

\subsection{Proof of Theorem \ref{thm:cartan}}

We divide the proof on several steps.

\subsubsection{}
First, we assume that the generalized root system is either finite or indefinite, as the affine case is discarded by \cite[Theorem 1.2]{AAH-rank2}. Moreover, we may assume that $\theta>2$, {\it ibid.}

\subsubsection{} Also, following \cite[Remark 3.2]{AAH-rank2}, we can assume that $\bq$ is of torsion class. More precisely, 
\begin{itemize}
	\item If $\bq$ is generic, then $\ba$ is of finite type by \cite{R quantum groups,AA}.
	\item If $\bq$ is semigeneric, then it is not of Cartan type. 
\end{itemize}
Let us develop the last assertion: If $\bq$ is semigeneric, then there is $i\in\I_\theta$ such that $q_{ii}\notin \G_\infty$ and $1\neq \widetilde{q}_{ij}\in\G_{\infty}$. But then $\widetilde{q}_{ij}\neq q_{ii}^a$, for any $a\in\Z_{\leq 0}$. 

\subsubsection{} Let $\ba$ be an indecomposable GCM of indefinite type. 

\begin{claim}
	We may assume that $\ba$ is compactly hyperbolic. 
\end{claim}
Indeed, if a given minor $\bm$ is of affine type and $\bq'$ is the submatrix of $\bq$ with the same entries as $\bm$, then $\bq'$ is of affine type $\bm$, so $\gdim\B_{\bq'}=\infty$. As $\B_{\bq'}$ is a subalgebra of $\B_{\bq}$, we have that $\gdim\B_{\bq}=\infty$ as well. The same argument allows us to consider only minimal matrices of indefinite type, that is with no minors of indefinite type. In other words, we can assume that $\ba$ is compactly hyperbolic. 
\medbreak

Now, if $\ba \in \Z^{\theta\times \theta}$ is compactly hyperbolic, then $\theta\leq 5$. The list of indecomposable GCM's which are compactly hyperbolic with rank $\theta>2$ contains exactly 31 matrices with $\theta=3$, three with $\theta=4$ and one with $\theta=5$. These are listed in \cite{C} and belong to
\begin{itemize}
	\item rows 1--31 of Tables 1 to 3 in loc. cit.~for $\theta=3$, 
	\item row 131 in Table 4 and rows 136, 137 in Table 5 for $\theta=4$,
	\item row 183 in Table 16 for $\theta=5$.
\end{itemize} 

\begin{claim}
	If $\ba$ is compactly hyperbolic, then $\dim \B_{\bq}=\infty$.
\end{claim}

This claim implies the theorem. We shall prove it in the series of Lemmas \ref{lem:cartan1}--\ref{lem:rank4} below.

\begin{notation}
	The enumeration of the rows in the tables in \cite{C} is continuous, (i.e.~Table $n+1$ starts at row $\ell+1$ if Table $n$ ends at row $\ell$); hence we will speak of a GCM $\ba$ which belongs to the corresponding row.
\end{notation}

Some of these matrices can be rapidly discarded, as follows.

\begin{lemma}\label{lem:cartan1}
	There is no matrix $\bq$ whose corresponding GCM is a generalized Cartan matrix $\ba$ of one of the following compactly hyperbolic types:
	\begin{enumerate}
		\item $\theta=3$, $\ba$ in rows $1,2,8,9,14,16,17,19,20,23$,
		\item $\theta=4$ and $\ba$ in row $131$,
		\item $\theta=5$.
	\end{enumerate}
\end{lemma}
\pf
In all of these cases, if we assume that there is $\bq=(q_{ij})$ with matrix $\ba$, we get that $q_{ij}=q_{ii}=1$, $i,j\in\I_\theta$; a contradiction. We work out a proof for $\ba$ in row 1, the rest follows similarly. In this case, the Dynkin diagram is
\[
\trianglerowone, \quad\text{hence}\quad \ba=\begin{pmatrix}
\,2&-1&-2\\-1&\,2&-1\\-1&-1&\,2
\end{pmatrix},
\]
which would give rise to a generalized Dynkin diagram of the following shape:
\[
\Dtriangle{q}{q^{-1}}{q}{q^{-1}=q^{-2}}{q^2}{q^{-2}}
\]
so $q=1$, a contradiction.
\epf

\begin{remark}\label{rem:AS}
	The same ideas in Lemma \ref{lem:cartan1} show that some diagrams of rank 3 in the list in \cite{C} can only come from a generalized Dynkin diagram with a root of 1 of a fixed order. The same holds for the rank 4 diagram on row 137, for which $\ord q=3$.
	
	These rank 3 matrices were spotted independently in \cite{AS-adv}. More precisely, the list, with corresponding order for $q$, is given in the following table:
	\begin{table}[H]
		\centering
		\begin{tabular}{ |c|c||c|c||c|c||c|c||c|c| } 
			\hline
			{\bf row} & $\ord q$ & {\bf row} & $\ord q$ &{\bf row} & $\ord q$ & {\bf row} & $\ord q$& {\bf row} & $\ord q$ \\
			\hline		\hline
			\bf5 & 3 & \bf6 & 5 & \bf7 & 5 & \bf12 & 8 & 1\bf3 & 7  		 \\ 
			\hline
			\bf15 & 11 & \bf18 & 4 & \bf21 & 17 & \bf22 & 7 & \bf24 & 26 	 \\ 
			\hline
		\end{tabular}
		\caption{Fixed order of root rows}
		\label{tableAS}
	\end{table}
\end{remark}

The following result will help us to deal with the content of Remark \ref{rem:AS}.

\begin{lemma}\label{lem:table-gral}
	Fix a braiding matrix $\bq=(q_{ij})_{i,j\in\I_3}$ such that
	\begin{align*}
	q_{jj}\widetilde{q}_{ij}&=1, & q_{kk}\widetilde{q}_{jk}&=q_{ii}\widetilde{q}_{ki}=1,
	\end{align*}
	for some choice of different $i,j,k\in\I_3$.
	\begin{enumerate}
		\item If $q_{ii}^2\widetilde{q}_{ij}=1$, then $\gdim\B_\bq=\infty$. 
		\item If $q_{ii}^3\widetilde{q}_{ij}=1$ and $q_{jj}\neq q_{kk}$, then $\gdim\B_\bq=\infty$. 
	\end{enumerate}
\end{lemma}
\pf
(1) We use Criterium 2 as in \S \ref{crit2}. Choose $\alpha=\alpha_1+\alpha_2+\alpha_3$, $\beta=2\alpha_i+\alpha_j$, $\omega=2\alpha_j-\alpha_i-\alpha_k$. Hence
\begin{align*}
q_{\alpha\beta}=\prod_{k\in\I_3}\widetilde{q}_{ki}^2\widetilde{q}_{kj}=q_{jj}\widetilde{q}_{ik}^2\widetilde{q}_{jk}=q_{jj}q_{ii}^{-2}q_{kk}^{-1}=q_{kk}^{-1}
\end{align*} 
and we get a diagram
\[
\Dchaintwo{1}{q_{kk}^{-1}}{q_{jj}},
\]
so $\gdim\B_\bq=\infty$. For (2), if  $\alpha=\alpha_1+\alpha_2+\alpha_3$, $\beta=3\alpha_i+\alpha_j$, so $\omega=3\alpha_j-\alpha_i-\alpha_k$, 
we get 
\[
\Dchaintwo{1}{q_{jj}q_{kk}^{-1}}{q_{jj}},
\]
hence $\gdim\B_\bq=\infty$.
\epf

\begin{lemma}\label{lem:table}
	If $\bq$ is such that either the corresponding GCM $\ba$ belongs to Table \ref{tableAS} in Remark \ref{rem:AS} or to row $137$, then $\gdim \B_\bq=\infty$.
\end{lemma}
\pf
If $\ba$ is not in rows $18,22$, then the result follows from items (1) and (2) of Lemma \ref{lem:table-gral}, 
with the right choice $ijk$ of indices $i,j,k\in\I_3$, as in the following table, which uses item (2) only for the last row:
\begin{table}[H]
	\centering
	\begin{tabular}{ |c|c||c|c||c|c||c|c| } 
		\hline
		{\bf row} & $ijk$ & {\bf row} & $ijk$ &{\bf row} & $ijk$ & {\bf row} & $ijk$ \\
		\hline		
		\bf5 & 321 & \bf6 & 321 & \bf7 & 231 & \bf12 & 321 	 \\ 
		\hline
		\bf13 & 123 & \bf15 & 321 &  \bf21 & 213 & \bf24 & 123 	 \\ 
		\hline
	\end{tabular}
\end{table}

Now, for rows $18$ and $22$ we  choose $\alpha=\alpha_1+\alpha_3$ and $\beta=\alpha_2+\alpha_3$ as in Criterium 3, \S \ref{crit3}, then we end up with the following diagrams:
\begin{align*}
& \Dchaintwo{q}{-1}{q^3}, \quad \ord q=4, &
& \Dchaintwo{q^{3}}{q^3}{q}, \quad \ord q=7;
\end{align*}
in both cases, $\gdim\B_\bq=\infty$.

For row $137$, we have $q=\zeta\in\G_3'$. We set $\alpha=\alpha_1$, $\beta=\alpha_2$ and $\gamma=\alpha_3+\alpha_4$ to obtain a triangle
\[
\tri{\zeta}{\zeta^2}{\zeta}{\zeta^2}{\zeta^2}{\zeta}\] 
For this, Criterium 1 for $\alpha=\alpha_1$ and $\beta=\alpha_2+2\alpha_3$ returns the segment
\[
\Dchaintwo{\zeta}{\zeta}{\zeta}
\]
which is not in the list, so $\gdim\B_{\bq}=\infty$.
\epf

\begin{lemma}\label{lem:rank3-triangle}
	If $\ba$ belongs to rows $3, 4, 10, 11$, then $\gdim\B_{\bq}=\infty$.
\end{lemma}
\pf
We proceed row by row.

{\it Row 3.} We use Criterium 1 as in \S \ref{crit1} for $\beta_1=\alpha_3$, $\beta_2=\alpha_1+\alpha_2$. We get
\[
\Dchaintwo{q^2}{q^{-4}}{q}
\]
which shows, by looking at \cite[Table 1]{H-rk3} that either $\ord q=3$ or 4 (disconnected dots). Now, if $\ord q=3$ and we choose 
$\beta_1=\alpha_1+\alpha_2+\alpha_3$, $\beta_2=2\alpha_2+\alpha_1$, as in Criterium 2, we get
\[
\Dchaintwo{q^{-1}}{q^{-4}}{q^2}\quad\equiv \quad \Dchaintwo{q^{-1}}{q^{-1}}{q^{-1}}
\]
which is not of finite type. On the other hand, if $\ord q=4$ and $\beta_1=2\alpha_2+\alpha_1$, $\beta_2=\alpha_1+\alpha_3$, we get
\[
\Dchaintwo{-1}{-q}{-q}
\]
and hence $\gdim\B_{\bq'}=\infty$. Here we have used the variant of Criterium 3, for $\omega=2\alpha_1-\alpha_2-2\alpha_3$, so $x_{13}$ and $x_{221}$ become primitive on $\K_{\omega}$, cf.~\eqref{eqn:del221}.

{\it Row 4.} 
By setting $\alpha=\alpha_3, \beta=\alpha_1+\alpha_2$ (Criterium 1) we obtain 
\[
\Dchaintwo{q}{q^{-3}}{q}
\]
so $\ord q=3,4$. We can discard $\ord q=4$ with $\alpha=2\alpha_1+\alpha_2, \beta=\alpha_1+\alpha_2+\alpha_3$, as in Criterium 2, and $\ord q=3$ with $\alpha=\alpha_1+\alpha_2, \beta=\alpha_3+\alpha_2$ (Criterium 3) as we obtain, in each case, the diagrams of affine type:
\begin{align*}
& \Dchaintwo{q^2}{q^{-2}}{q^{-1}}, & & \Dchaintwo{q^{-1}}{q^{-1}}{q^{-1}}.
\end{align*}

{\it Row 10.} If $\alpha=\alpha_1+\alpha_3$, $\beta=\alpha_2+\alpha_3$ (Criterium 3), then we obtain
\[
\Dchaintwo{q}{q}{q},
\]
As $\ord q\geq 4$, we have that $\gdim\B_\bq=\infty$.

{\it Row 11.} Again, $\ord q> 3$ and the choice $\alpha=\alpha_1+\alpha_3$, $\beta=\alpha_2$ (Criterium 1) gives $\ord q=6$. 
But $\ord q\neq 6$ by considering $\alpha=\alpha_2$, $\beta=2\alpha_1+\alpha_3$.
\epf 

Next we turn our attention to the remaining diagrams of rank 3, which are not triangles.

\begin{lemma}\label{lem:rank3-not-triangle}
	If $\ba$ belongs to rows $25, \dots, 31$, then $\gdim\B_{\bq}=\infty$.
\end{lemma}
\pf
Matrices $\bq$ with $\ba$ in rows 28, respectively 31, are
\begin{align*}
& \Dchainthree{q^2}{q^{-6}}{q^6}{q^{-6}}{q^3}, && \ord q\ne 2,3,4,6,
\\
&\Dchainthree{q}{q^{-3}}{q^3}{q^{-3}}{r}, && \ord q,\ord r \ne 2,3, q^3=r^3.
\end{align*}
For $\alpha=2\alpha_1+\alpha_2$, $\beta=\alpha_2+\alpha_3$ we get the diagrams
\begin{align*}
& \Dchaintwo{q^2}{q^{-6}}{q^3},
&
& \Dchaintwo{q}{q^{-3}}{r},
\end{align*}
both of indefinite Cartan type, so $\gdim\B_{\bq}=\infty$. Rows $25$, $26$, $27$, $29$ and $30$ are discarded by setting 
$\alpha=\alpha_1+\alpha_2$, $\beta=\alpha_2+\alpha_3$, using Criterium 3. Indeed the corresponding diagrams are
\begin{align*}
& \Dchaintwo{q^2}{q^{-4}}{q}, & 
& \Dchaintwo{q}{q^{-3}}{q^6}, & 
& \Dchaintwo{q}{q^{-3}}{q}, & 
& \Dchaintwo{q^3}{q^{-6}}{q}, & 
& \Dchaintwo{q}{q^{-4}}{q}.
\end{align*}
Due to the restrictions on $\ord q$ for each case, these diagrams are connected, and moreover of indefinite Cartan type.
\epf

We end up this analysis with the only diagram left, of rank 4.
\begin{lemma}\label{lem:rank4}
	If $\ba$ is the GCM of row $136$, then $\gdim\B_\bq=\infty$.
\end{lemma}
\pf
Take $\beta_1=\alpha_1$, $\beta_2=\alpha_2+\alpha_3$ and $\beta_3=\alpha_4$. We obtain the rank 3 diagram associated to row $3$, and the statement follows by Lemma \ref{lem:rank3-triangle}.
\epf

This finishes the proof of the claim, hence of Theorem \ref{thm:cartan}.

\section{Rank three case}\label{sec:rank3}

In this section, we give a positive answer to Conjecture \ref{conjecture} for $\theta=3$. 

To this end, we will collect a series of lemmas, which are stated after the proof for the sake of readability, that deal with diagrams $d$ constructed by the pasting of two (a \qt{line}) or three (a \qt{triangle}) diagrams of rank two with finite root system.  

We say that a diagram $d$ of rank two in the list is \qt{finite} if all of its entries are concrete roots of 1 in $\G_\infty$
(that is, $d$ belongs to one of the rows 6--10 or 12--17) and call it \qt{parametric} otherwise. Parametric diagrams are indeed in the one-parameter families mentioned in \cite{AAY}. 

The main result of this section reads as follows.
\begin{theorem}\label{thm:rank3}
	Let $(V,c^{\bq})$ be a braided vector space of diagonal type with $\dim V=3$. If $\gdim\B_{\bq}<\infty$, then the root system of $\bq$ is
	finite. 
\end{theorem}
\pf
Assume that $\gdim\B_{\bq}<\infty$. 
The diagram $d$ associated to $\bq$ is either a line or a triangle; in any case, each subdiagram of rank 2 belongs to the list by \cite[Theorem 4.1]{AAH-rank2}. In particular, these subdiagrams can be finite or parametric.

Assume first that $d$ is a line. Then, it is obtained by gluing either two parametric diagrams of rank 2, a parametric diagram with a finite one, or two finite ones. 
By Lemmas \ref{lem:line-1param-1finite} and \ref{lem:line-two-parametric} we have that $\bq$ is in the list or the parametric subdiagrams (if any) are evaluated in the finite set $\Gf$ from \eqref{eq:defn-G-f}. But in this last case, Lemma \ref{lem:gaplines} shows that $\bq$ is in the list as well.

When $d$ is a triangle, the situation is analogous. Indeed, first we use Lemmas \ref{lem:triang-2param-1finite}, \ref{lem:triang-1param-2finite} and \ref{lem:triang-3param} to deduce that $\bq$ is in the list or the parametric subdiagrams (if any) are evaluated again in $\Gf$. Then Lemma \ref{lem:gaptriangles} implies that $\bq$ is in the list.
\epf

\subsection{On parametric subdiagrams}\label{subsec:eval-param}

Next, we deal with diagrams of rank 3 with at least one parametric subdiagram. These subdiagrams can be evaluated in a root of 1 or not.
In the first case, we will show that it is enough to test the results when evaluating on the set
\begin{align}\label{eq:defn-G-f}
\Gf &:= \G_{14}\cup \G_{18} \cup \G_{20}\cup \G_{24} \cup \G_{30}.
\end{align}
This reduction is done in the following five lemmas, which deal with all combinations of finite and parametric subdiagrams of rank two.
That is, we have, respectively, the following cases
\begin{enumerate}
	\item A line with a parametric subdiagram and a finite subdiagram;
	\item A line with two parametric subdiagrams;
	\item A triangle with one parametric subdiagram and two finite ones; 
	\item A triangle with two parametric subdiagrams and a finite one;
	\item A triangle with three parametric subdiagrams.
\end{enumerate}

\begin{remark}
	The entries of each finite diagram in the list belong to $\Gf$.
\end{remark}

After this reduction, we are left with a finite (though large) collection of diagrams with labels in $\Gf$, that we attack with \gap\, in \S \ref{subsec:gap}.

We set the following notation for a \emph{parametric} diagram; i.e.~a function $\mathtt{G}:\C^{\times}\to (\C^{\times})^3$, where $\mathtt{G}(q)$ is one of the following functions:
\begin{align*}
\Ados{q}&=\Dchaintwo{q}{q^{-1}}{q}; &
\Bdos{q}&=\Dchaintwo{q}{q^{-2}}{q^2}; &
\Bdost{q}&=\Dchaintwo{q^2}{q^{-2}}{q}; 
\\
\Gdos{q}&=\Dchaintwo{q}{q^{-3}}{q^3}; &
\Gdost{q}&=\Dchaintwo{q^3}{q^{-3}}{q}; &
\sAuno{q}&=\Dchaintwo{-1}{q}{-1};
\\
\sAdos{q}&=\Dchaintwo{q}{q^{-1}}{-1}; &
\sAdost{q}&=\Dchaintwo{-1}{q^{-1}}{q}; &
\sB{q}&=\Dchaintwo{q}{q^{-2}}{-1}; &
\\
\sBt{q}&=\Dchaintwo{-1}{q^{-2}}{q}; &
\Bst{q}&=\Dchaintwo{q}{q^{-1}}{\zeta}; &
\Bstt{q}&=\Dchaintwo{\zeta}{q^{-1}}{q}.
\end{align*}

\begin{lemma}\label{lem:line-1param-1finite}
	Let $\bq$ be a line such that $\bq_{|\I_2}$ is parametric of type $\mathtt{G}$ and $\bq_{|\I_{2,3}}$ is finite. Then either $\bq_{|\I_2}=\Gt{q}$ for some $q\in\Gf$ or else $\gdim \B_{\bq}=\infty$.
\end{lemma}
\pf
As $\bq$ is made by glueing $\Gt{q}$ with a finite diagram $\Ft=\Dchaintwo{\Ft_1}{\Ft_2}{\Ft_3}$, the last entry $\Gt{q}_3$ of $\Gt{q}$ is $\Ft_1$, and $\Ft_i \in \Gf$ for all $i\in\I_3$.

If $\Gt{q}_3=q$, then the claim follows. Now we check the remaining cases:
\smallbreak

\begin{description}[leftmargin=*]
	\item[$\Bdos{q}$] Applying \one\,  to $\alpha_1+\alpha_2$, $\alpha_3$ we obtain the diagram $d=\Dchaintwo{q}{\Ft_2}{\Ft_3}$. Suppose that $\gdim \B_{\bq}<\infty$: If $d$ is finite, then $q\in\Gf$. If $q\Ft_2=1$ or $q=\zeta\in\G_3'$, then the same holds. Otherwise, either $q^2\Ft_2=1$ with $\Ft_3\in\{-1,\Ft_2^{-1}\}$ or $q^3\Ft_2=1=\Ft_2\Ft_3$. By inspection, a finite diagram such that $\Ft_3\in\{-1,\Ft_2^{-1}\}$ satisfy $\Ft_1\in\G_{6}\cup\G_{10}$, in which case $q\in\Gf$ since $\Ft_1=q^2$, or $\Ft$ is one of the following diagrams:
	\begin{align*}
	& \Dchaintwo{\zeta^5}{\zeta^9}{-1}, \, \zeta\in\G_{12}'; & 
	& \Dchaintwo{-\zeta^2}{\zeta}{-1}, \, \zeta\in\G_{9}';
	\\
	& \Dchaintwo{\zeta}{\zeta^{-5}}{-1}, \, \zeta\in\G_{24}'; & 
	& \Dchaintwo{\zeta}{\zeta^{-3}}{-1}, \, \zeta\in\G_{20}'.
	\end{align*}
	Applying \one\,  to $\alpha_1$, $\alpha_2+\alpha_3$ we obtain the diagram $\Dchaintwo{q}{q^{-2}}{-\Ft_1\Ft_2}$, which does not belong to Heckenberger's list.
	
	\smallbreak
	
	\item[$\Gdos{q}$] As before, \one\,  for $\alpha_1+\alpha_2$, $\alpha_3$ gives the diagram $\Dchaintwo{q}{\Ft_2}{\Ft_3}$ and the proof follows as for $\Bdos{q}$
	\smallbreak
	
	\item[$\sAuno{q}$] Here $\Ft_1=-1$, so $\Ft_3\ne -1$. Applying \one\,  to $\alpha_1+\alpha_2$, $\alpha_3$ we obtain the diagram $d=\Dchaintwo{q}{\Ft_2}{\Ft_3}$. Suppose that $\gdim \B_{\bq}<\infty$: If $d$ is finite, then $q\in\Gf$, as well as if $q\Ft_2=1$ or $q=\zeta\in\G_3'$. Otherwise, $1=\Ft_3\Ft_2$, so $\Ft$ is parametric, a contradiction.
	
	\smallbreak
	
	\item[$\sAdos{q}$] Applying \one\,  to $\alpha_1$, $\alpha_2+\alpha_3$ we get the diagram $d=\Dchaintwo{q}{q^{-1}}{-\Ft_2\Ft_3}$. Suppose that $\gdim \B_{\bq}<\infty$: If $d$ is finite, then $q\in\Gf$. As $\Ft$ is not parametric, $\Ft_2\Ft_3\ne 1$; the remaining possibilities are $q\in \{-\Ft_2\Ft_3, \Ft_2^2\Ft_3^2, -\Ft_2^3\Ft_3^3\}\subset \Gf$, or $-\Ft_2\Ft_3\in\G_3'$. By inspection, $\Ft$ is one of the following diagrams:
	\begin{align*}
	& \Dchaintwo{-1}{-\zeta}{\zeta}, \, \zeta\in\G_{3}'; & 
	& \Dchaintwo{-1}{\zeta}{-\zeta^3}, \, \zeta\in\G_{12}'; & 
	& \Dchaintwo{-1}{-\zeta^3}{-\zeta^{-1}}, \, \zeta\in\G_{12}'; 
	\\
	& \Dchaintwo{-1}{\zeta}{-\zeta^2}, \, \zeta\in\G_{9}'; & 
	& \Dchaintwo{-1}{\zeta^{-5}}{\zeta}, \, \zeta\in\G_{24}'; & 
	& \Dchaintwo{-1}{-\zeta^2}{\zeta^3}, \, \zeta\in\G_{15}'.
	\end{align*}
	Applying \one\,  to $\alpha_1$, $\alpha_2+2\alpha_3$ to the last five diagrams we get the diagram $D'=\Dchaintwo{q}{q^{-1}}{-\Ft_2^2\Ft_3^4}$: either $D'$ is not in the list (so $\gdim \B_{\bq}=\infty$), $D'$ is finite (in which case $q\in\Gf$), or $q\in \{-\Ft_2^2\Ft^4, \Ft_2^4\Ft^8, -\Ft_2^6\Ft^9\}\subset \Gf$. For the first diagram, we reflect on vertex 2 and obtain a triangle with the subdiagram
	$\Dchaintwo{-1}{-\zeta q^{-1}}{\zeta^2}$. 
	By scanning the list, we see that this is not the list unless $q\in \G_{18}\cup \G_{24}\cup \G_{30}\subset \Gf$.
	
	\smallbreak
	
	\item[$\sB{q}$] Applying \one\,  to $\alpha_1+\alpha_2$, $\alpha_3$ we have the diagram $\Dchaintwo{-q^{-1}}{\Ft_2}{\Ft_3}$ and the proof follows as for $\Bdos{q}$.
	
	\smallbreak
	
	\item[$\Bst{q}$] In this case, \one\,  with $\alpha_1$, $\alpha_2+\alpha_3$ gives $d=\Dchaintwo{q}{q^{-1}}{\zeta\Ft_2\Ft_3}$. Now, if $d$ is finite, then $q\in\Gf$. If $d$ is parametric, then either $q=(\zeta\Ft_2\Ft_3)^m$, $m\leq 3$, in which case $q\in\Gf$, or either $\zeta\Ft_2\Ft_3\in\{-1,\zeta^2\}$, that is $\Ft=\Dchaintwo{\zeta}{-\zeta}{-1}$. We apply  \one\,  with $\alpha_1+2\alpha_2$, $\alpha_3$ and obtain  
	\begin{align*}
	& \Dchaintwo{\zeta q^{-1}}{\zeta^2}{-1} 
	\end{align*}
	By inspection, this diagram cannot be finite; hence it is parametric which forces $q\in\G_6$ or it is not in Heckenberger's list, giving $\gdim\B_{\bq}=\infty$. 
\end{description}

This shows the result.
\epf

\begin{lemma}\label{lem:line-two-parametric}
	Let $\bq$ be a line such that $\bq_{|\I_2}$ and $\bq_{|\I_{2,3}}$ are parametric of type $\mathtt{G}$, $\mathtt{G}'$. Then either $\bq$ belongs to Heckenberger's list, or $\bq_{|\I_2}=\Gt{q}$, $\bq_{|\I_{2,3}}=\Gtp{r}$ for some $q,r\in\Gf$ or else $\gdim \B_{\bq}=\infty$.
\end{lemma}
\pf
If $\bq$ is of Cartan type, then either $\bq$ belongs to Heckenberger list or or else $\gdim \B_{\bq}=\infty$ by Theorem \ref{thm:cartan}. Hence we may assume that $\bq_{|\I_2}$ is not of Cartan type. We 
prove the statement case by case depending on the type of $\bq_{|\I_2}$. At each step we will not consider the possibility that $\bq_{|\I_{2,3}}$ has a diagram already analyzed: unless the vertices 1 and 3 are exchanged, $\bq$ has already been considered.

\medspace

\begin{description}[leftmargin=*]
	\item[$\diamondsuit\sAdos{q}$] The possible $\mathtt{G}'$ are 
	\begin{align*}
	\Ados{-1}&, &
	\Bdost{\xi}&, \xi\in\G_4', &
	\Gdost{\xi}&, \xi\in\G_6', &
	\sAuno{r}&, &
	\sAdost{r}&, &
	\sBt{r}&, &
	\Bst{-1}&.
	\end{align*}
	For $\Ados{-1}$ and $\sAdost{r}$, $\bq$ belongs to rows 9, 10 or 11, depending on $r$.
	
	If $\mathtt{G}'$ is either $\Bdost{\xi}$, $\xi\in\G_4'$, or $\Gdost{\xi}$, $\xi\in\G_6'$, then \one\,  for $\alpha_1$, $\alpha_2+\alpha_3$ gives the diagram $\Dchaintwo{q}{q^{-1}}{\xi}$: either this diagram is not in the list (so $\gdim \B_{\bq}=\infty$), or it is finite (which implies $q \in\Gf$), or parametric, in which case $q\in\{\xi, \xi^2, \xi^3\}$. Similar for $\Bst{-1}$, since \one\,  for $\alpha_1$, $\alpha_2+\alpha_3$ gives the diagram $\Dchaintwo{q}{q^{-1}}{-\zeta}$.
	
	For $\bq_{|\I_{2,3}}=\sAuno{r}$, \one\,  for $\alpha_1$, $\alpha_2+\alpha_3$ gives the diagram $\Dchaintwo{q}{q^{-1}}{r}$. Assume that this diagram is in the list: either $q,r\in\Gf$ if the diagram is finite, or $q\in \{r,r^2,r^3\}$, or $r=\zeta$. The cases $q=r,r^2,r^3$ belong, respectively, to rows 4, 6 and 7. For $r=\zeta$, the diagram of $\rho_2\bq$ is
	$$ \tri{-1}{q}{-1}{\zeta^2}{\zeta}{\zeta q^{-1}} .$$
	Looking at the subdiagram of 1 and 3, we have three possibilities: $q\in\Gf$ if this diagram is finite, $q\in\G_6$ if it is parametric, or $\gdim \B_{\bq}=\infty$.
	
	For $\bq_{|\I_{2,3}}=\sBt{r}$, \one\,  for $\alpha_1$, $\alpha_2+\alpha_3$ gives $\Dchaintwo{q}{q^{-1}}{-r^{-1}}$. Assume that this diagram is in the list: either $q,r\in\Gf$ if the diagram is finite, or $q\in \{-r^{-1},r^{-2},-r^{-3}\}$, or $r=-\zeta^2$.
	The cases $q=r^2,-r^3$ belong, respectively, to rows 5 and 7. For $q=-r^{-1}$, the subdiagram of $\rho_2\bq$ with vertices 1,3 is $\Dchaintwo{-1}{-r^{-1}}{-r^{-1}}$, which either is not in the list, or $r\in\G_{12}$, in which case $q\in\G_{24}$. For $r=-\zeta^2$, the subdiagram of $\rho_2\bq$ with vertices 1,3 is $\Dchaintwo{-1}{q^{-1}\zeta^2}{\zeta}$, which either is not in the list, or $q\in\G_{12}$.
	
	\medspace
	
	\item[$\diamondsuit\sAdost{q}$] The possible $\mathtt{G}'$ are 
	\begin{align*}
	\Ados{q}&, &
	\Bdos{q}&, &
	\Bdost{r}&, \, q=r^2,&
	\sB{q}&, &
	\Bstt{r}&, q=\zeta,
	\\
	\sAdos{q}&, &
	\Gdos{q}&, &
	\Gdost{r}&, \, q=r^3,&
	\Bst{q}&. & &
	\end{align*}
	The first, the second, the third, the sixth and the seventh diagrams appear, respectively, in rows 4, 6, 5, 8 and 7. 
	
	For $\sB{q}$, \one\,  for $\alpha_1$, $\alpha_2+\alpha_3$ gives $\Dchaintwo{-1}{q^{-1}}{-q^{-1}}$. If this diagram is finite, then $q\in\Gf$. If it is parametric, then either $-q^{-2}=1$ or $q^{-3}=1$. Similarly, for $\Bst{q}$, \one\,  for $\alpha_1$, $\alpha_2+\alpha_3$ gives the diagram $\Dchaintwo{-1}{q^{-1}}{\zeta}$, and again $q\in\Gf$.
	For $\Bstt{r}$, $q=\zeta$, the subdiagram of $\rho_2\bq$ with vertices 1,3 is $\lin{-1}{r^{-1}}{r^{-1}\zeta}$: either it is not in the list or $r\in\Gf$.
	Finally, for $\Gdost{r}$, $q=r^3$, \one\,  for $\alpha_1$, $\alpha_2+\alpha_3$ gives the diagram $\Dchaintwo{-1}{r^{-3}}{r}$: this diagram is not parametric, so either $\gdim \B_{\bq}=\infty$ or $r\in\Gf$.
	
	\medspace
	
	\item[$\diamondsuit\sAuno{q}$] The possible $\mathtt{G}'\ne \sAdos{r}, \sAdost{r}$ are 
	\begin{align*}
	\Ados{-1}&, &
	\sAuno{r}&, &
	\sB{r}&, &
	\Bdost{\xi}&, \xi\in\G_4', &
	\Gdost{\xi}&, \xi\in\G_6',&
	\Bstt{-1}&.
	\end{align*}
	For $\Ados{-1}$, \one\,  for $\alpha_1+\alpha_2$, $\alpha_3$ gives the diagram $\Dchaintwo{q}{-1}{-1}$: either $\ord q\in\{3,4,6\}$ or this diagram does not belong to the list.
	
	For $\sAuno{r}$, \one\,  for $\alpha_1+\alpha_2$, $\alpha_3$ and for $\alpha_1$, $\alpha_2+\alpha_3$ give, respectively, the diagrams $\Dchaintwo{q}{r}{-1}$ and $\Dchaintwo{-1}{q}{r}$: either $q,r\in\Gf$, or $qr=1$ (which appears in row 8), or $q^2r=1=qr^2$ (so $q=r\in\G_3'$, and $\bq$ appears in Row 15), or $\gdim\B_{\bq}=\infty$.
	
	For $\sB{r}$, \one\,  for $\alpha_1$, $\alpha_2+\alpha_3$ gives the diagram $\Dchaintwo{q}{r^{-2}}{r}$: either $q,r\in\Gf$, or $q=r^2$ (which appears in row 5), or $\gdim\B_{\bq}=\infty$.
	
	For $\Bdost{\xi}$, $\xi\in\G_4'$, \one\,  for $\alpha_1+\alpha_2$, $\alpha_3$ gives $\Dchaintwo{q}{-1}{\xi}$, which does not appear in the list. The same happens for $\Gdost{\xi}$, $\xi\in\G_6'$, and $\Bstt{-1}$.
	
	\medspace
	
	\item[$\diamondsuit\sB{q}$] There are 5 possible $\mathtt{G}'\ne \sAdos{r}, \sAdost{r}, \sAuno{r}$. One is $\Ados{-1}$: \one\,  for $\alpha_1+\alpha_2$, $\alpha_3$ gives the diagram $\Dchaintwo{-q^{-1}}{-1}{-1}$, so either $\ord q\in\{3,4,6\}$ or this diagram does not belong to the list. 
	
	For $\Bdost{\xi}$, $\xi\in\G_4'$, \one\,  for $\alpha_1+\alpha_2$, $\alpha_3$ gives $\Dchaintwo{-q^{-1}}{-1}{\xi}$, which does not appear in the list. The same with $\Gdost{\xi}$, $\xi\in\G_6'$, and $\Bstt{-1}$. Finally, for $\sBt{r}$, \one\,  for $\alpha_1+\alpha_2$, $\alpha_3$ and for $\alpha_1$, $\alpha_2+\alpha_3$ give, respectively, the diagrams $\Dchaintwo{-q^{-1}}{r^{-2}}{r}$ and $\Dchaintwo{q}{q^{-2}}{-r^{-1}}$: either $q,r\in\Gf$, or $q^2r=qr^2=-1$ (so $q=r\in\G_6'$), or $\gdim\B_{\bq}=\infty$.
	
	\medspace
	
	\item[$\diamondsuit\sBt{q}$] The possible $\mathtt{G}'\ne \sAdos{r}, \sAdost{r}$ are 
	\begin{align*}
	&\Ados{q}, &
	&\sB{q}, &
	&\Bdos{q}, &
	&\Bdost{r}, q=r^2,
	\\
	&\Gdos{q}, &
	&\Gdost{r}, q=r^3,&
	&\Bstt{r}, q=\zeta, &
	&\Bst{q}.
	\end{align*}
	\one\,  for $\alpha_1+\alpha_2$, $\alpha_3$ gives, respectively, the diagrams 
	\begin{align*}
	& \Dchaintwo{-q^{-1}}{q^{-1}}{q}, &
	& \Dchaintwo{-q^{-1}}{q^{-2}}{-1}, &
	& \Dchaintwo{-q^{-1}}{q^{-2}}{q^2}, &
	& \Dchaintwo{-r^{-2}}{r^{-2}}{r},
	\\
	& \Dchaintwo{-q^{-1}}{q^{-3}}{q^3}, &
	& \Dchaintwo{-r^{-3}}{r^{-3}}{r}, &
	& \Dchaintwo{-\zeta^{2}}{r^{-1}}{r}, &
	& \Dchaintwo{-q^{-1}}{q^{-1}}{\zeta}.
	\end{align*}
	For each one of these diagrams $d=\Dchaintwo{d_1}{d_2}{d_3}$ we perform the following tripartite analysis:
	\begin{itemize}[leftmargin=*]
		\item If the diagram is finite, then $q,r\in\Gf$. Indeed, $q\in \{r,r^2,r^3\}$ or $q\in\G_3'$: in any case $r\in\Gf$, thus $q\in\Gf$ too.
		\item If it is parametric, then we see that $d_1\ne -1$ and therefore there is $m\in\I_3$ such that $d_1^md_2=1$. This automatically implies that $q, r\in\Gf$.
		\item Otherwise, $\gdim \B_{\bq}=\infty$.
	\end{itemize}

	\medspace
	
	\item[$\diamondsuit\Bst{q}$] The possible $\mathtt{G}'\ne \sAdos{r}, \sAdost{r}, \sB{r},\sBt{r}$ are 
	\begin{align*}
	&\Ados{\zeta}, &
	&\Bdos{\zeta}, &
	&\Bdost{\pm \zeta^2}, &
	&\Gdost{\xi}, \xi^3=\zeta,&
	&\Bstt{r}.
	\end{align*}
	\one\,  for $\alpha_1+2\alpha_2$, $\alpha_3$ gives, respectively:
	\begin{align*}
	& \Dchaintwo{\zeta q^{-1}}{\zeta}{\zeta}, &
	& \Dchaintwo{\zeta q^{-1}}{\zeta^2}{\zeta^2}, &
	& \Dchaintwo{\zeta q^{-1}}{\zeta}{\pm \zeta}, &
	& \Dchaintwo{\zeta q^{-1}}{\xi^3}{\xi}, &
	& \Dchaintwo{\zeta q^{-1}}{r^{-2}}{r}.
	\end{align*}
	The first four diagrams are discarded with the same tripartite analysis as above. For the last, we have that it is in the list if and only if either $\zeta q^{-1}=r^2$, or $\zeta q^{-1}=-1$, or $\zeta q^{-1}=\xi^3$, $r=-\xi$ for some $\xi\in\G_9'$. Analogously, \one\,  for $\alpha_1$, $2\alpha_2+\alpha_3$ gives $\Dchaintwo{q}{q^{-2}}{\zeta r^{-1}}$, which is in the list if and only if either $\zeta r^{-1}=q^2$, or $\zeta r^{-1}=-1$, or $\zeta r^{-1}=\eta^3$, $r=-\eta$ for some $\eta\in\G_9'$.
	The combination of these possibilities gives $q,r\in\Gf$.
	\medspace
	
	\item[$\diamondsuit\Bstt{q}$] The possible $\mathtt{G}'\ne \sAdos{r}, \sAdost{r}, \sB{r},\sBt{r}$ are 
	\begin{align*}
	&\Ados{q}, &
	&\Bdos{q}, &
	&\Bdost{r}, q=r^2 &
	&\Gdos{q}, &
	&\Gdost{r}, q=r^3,&
	&\Bst{q}.
	\end{align*}
	On the one hand, for $\Ados{q}$, $\Bdos{q}$ and $\Gdos{q}$ we apply \one\,  for $2\alpha_1+\alpha_2$, $\alpha_3$ and obtain
	$\Dchaintwo{\zeta q^{-1}}{q^{-m}}{q^m}$ for $m=1,2,3$, respectively.
	On the other, for $\Bdost{r}$ with $q=r^2$, $\Gdost{r}$ with $q=r^3$ and $\Bst{q}$,
	\one\,  for $\alpha_1$, $\alpha_2+\alpha_3$ gives, respectively:
	\begin{align*}
	& \Dchaintwo{\zeta}{r^{-2}}{r}, &
	& \Dchaintwo{\zeta}{r^{-3}}{r}, &
	& \Dchaintwo{\zeta}{q^{-1}}{\zeta}.
	\end{align*}
	Again, we apply the tripartite analysis.
\end{description}

This ends the proof.\epf

\begin{lemma}\label{lem:triang-2param-1finite}
	Let $\bq$ be a triangle such that $\bq_{|\I_2}$ is finite and both $\bq_{|\I_{2,3}}$ and  $\bq_{|\{1,3\}}$ are parametric. Then $\bq_{|\I_2}=\Gt{q}$ and $\bq_{|\{1,3\}}=\Gtp{r}$ for some $q,r\in\Gf$ or else $\gdim \B_{\bq}=\infty$.
\end{lemma}
\pf
We fix $\zeta\in\G_3'$. We deal with a triangle of the form
\[
\tri{\Gt{r}'_3=\Ft_1}{\Ft_2}{\Ft_3=\Gt{q}_1}{\Gt{q}_2}{\Gt{q}_3=\Gt{r}'_1}{\Gt{r}'_2}
\]
If $\Gt{q}_1=q$ and $\Gt{r}'_3=r$, then $q,r\in\Gf$. 
If $\Gt{q}_1=q$ and $\Gt{r}'_1=r$, then $q \in \Gf$ and $r\in\{-1,\zeta^{\pm1},q,q^2,q^3\}$, 
which implies that $r\in\Gf$.
The same holds if, conversely, $\Gt{q}_3=q$ and $\Gt{r}'_3=r$.

We are left with $q=\Gt{q}_3=\Gt{r}'_1=r$. Hence $\Ft_1,\Ft_3\in\{-1,\zeta^{\pm1},q,q^2,q^3\}$. We analyze this case by case.
\begin{enumerate}[leftmargin=*]
	\item If either $\Ft_1=q$ or $\Ft_3=q$, then $q\in\Gf$.
	\item If $\Ft_1=\Ft_3$ then by inspection we have that  $\Ft_2=\beta\in\G_{12}'$ and $\Ft_1=\Ft_3=\beta^8\in \G_{3}'$ (diagram \tt{d9.1}). Hence the triangle is
	\[
	\tri{\beta^8}{\beta}{\beta^8}{q^{-1}}{q}{q^{-1}}
	\]
	Applying \one\,  to $\alpha_1+\alpha_2$ and $\alpha_3$ we obtain the line 
	\begin{align*}
	& \Dchaintwo{\beta^5}{q^{-2}}{q} 
	\end{align*}
	If this diagram is finite, then $q\in\Gf$; if it is parametric, then it is $\Bdost{q}$, for $q^2=\beta^5\in\G'_{12}$ and thus $q\in\Gf$. Otherwise, it is not in the list.
	\item Assume $\Ft_1=q^2$, $\Ft_3=q^3$, equivalently $\Ft_3^2=\Ft_2^3$.  By inspection of the list, we see that $q\in\Gf$.
	
	\item When $\Ft_1=q^2$ and $\Ft_3=\zeta$, we obtain the triangles:
	\[
	\tri{q^2}{\Ft_2}{\zeta}{q^{-1}}{q}{q^{-2}}
	\]
	We have that either $q^2\in\G_{12}\cup\G_{15}$, so $q\in\Gf$, or $\Ft=\Dchaintwo{-\beta}{\beta^{-2}}{\beta^3}$, $\beta\in\G_{9}'$ (\tt{d10.1}) or $\Ft=\Dchaintwo{-\beta}{-\beta^{-3}}{\beta^5}$, $\beta\in\G_{15}'$ (\tt{d16.1}). By renaming the variables (picking the right $q$), we have two possible triangles, namely,
	\begin{align*}
	&\tri{\xi^2}{\xi^{-4}}{-\xi^6=\zeta}{\xi^{-1}}{\xi}{\xi^{-2}}, \xi\in\G_{36}' & \tri{\xi^2}{\xi^{-6}}{-\xi^{10}=\zeta}{\xi^{-1}}{\xi}{\xi^{-2}}, \xi\in\G_{60}'
	\end{align*} 
	Applying \one\,  for $\alpha_1$ and $\alpha_2+\alpha_3$ in both cases we obtain the lines
	\begin{align*}
	&   \Dchaintwo{\xi^2}{\xi^{-6}}{-\xi^{6}}, \xi\in\G_{36}' & \Dchaintwo{\xi^{2}}{\xi^{-8}}{\xi^{-20}}, \xi\in\G_{60}';
	\end{align*}
	which are not in the list.
	
	\item Now we assume $\Ft_1=q^2$, $\Ft_3=-1$, that is
	\[
	\tri{q^2=\Ft_1}{\Ft_2}{-1}{q^{-m}\qquad m=1,2}{q}{q^{-2}}
	\] 
	A new inspection gives either $q\in\Gf$, $\Ft_1=-\beta^2$, $\Ft_2=\beta$, $\beta\in \G'_9$, (\tt{d10.3}) or else $\Ft_1=\beta$, $\Ft_2=\beta^{-3}$, $\beta\in \G'_{20}$, (\tt{d15.1}). Hence we have the following possibilities:
	\begin{align*}
	&\tri{\xi^2}{\xi^{-8}}{-1}{\xi^{-1}}{\xi}{\xi^{-2}},  & \tri{\xi^2}{\xi^{-8}}{-1}{\xi^{-2}}{\xi}{\xi^{-2}}, \xi\in\G_{36}';\\
	&\tri{\xi^2}{\xi^{-6}}{-1}{\xi^{-1}}{\xi}{\xi^{-2}},   & \tri{\xi^2}{\xi^{-6}}{-1}{\xi^{-2}}{\xi}{\xi^{-2}}, \xi\in\G_{40}'.
	\end{align*} 
	For these triangles, we apply \one\,  for $\alpha_1$, $\alpha_2+\alpha_3$ and obtain the lines
	\begin{align*}
	&   \Dchaintwo{\xi^{2}}{-\xi^{8}}{-1}, && \Dchaintwo{\xi^{2}}{-\xi^{8}}{-\xi^{-1}}, & \xi\in\G_{36}';\\
	&   \Dchaintwo{\xi^{2}}{-\xi^{12}}{-1}, && \Dchaintwo{\xi^{2}}{-\xi^{12}}{-\xi^{-1}}, & \xi\in\G_{40}';
	\end{align*}
	which are not in the list. 
	
	\item In this case, $\Ft_1=q^3$, $\Ft_3=\zeta$, so
	\[
	\tri{q^3=\Ft_1}{\Ft_2}{\zeta}{q^{-1}}{q}{q^{-3}}
	\]
	Then, by inspection, either $q\in\Gf$ or $\Ft$ is of types \tt{d10.1} or \tt{d16.1} as above and we end up with the following triangles:
	\begin{align*}
	&\tri{\xi^3}{\xi^{-6}}{-\xi^9}{\xi^{-1}}{\xi}{\xi^{-3}}, \xi\in\G'_{54}, & \tri{\xi^3}{\xi^{-9}}{-\xi^{15}}{\xi^{-1}}{\xi}{\xi^{-3}}, \xi\in\G'_{90},
	\end{align*}
	Again, \one\,  for $\alpha_1$, $\alpha_2+\alpha_3$ gives lines which are not in the list, namely, 
	\begin{align*}
	&   \Dchaintwo{\xi^{3}}{-\xi^{8}}{-\xi^{9}}, \xi\in\G'_{54}  & \Dchaintwo{\xi^{3}}{-\xi^{33}}{-\xi^{15}},  \xi\in\G_{90}'.
	\end{align*}
	
	\item We have $\Ft_1=q^3$, $\Ft_3=-1$. Then either $q \in\Gf$ or $\Ft$ is one of \tt{d9.3}, \tt{d10.3}, \tt{d15.1}, \tt{d17.1} or \tt{d17.2}. This gives the following possibilities:
	\begin{align*}
	&\tri{\xi^3}{\xi^{-9}}{-1}{\xi^{-1}}{\xi}{\xi^{-3}}, \xi\in\G_{36}', & \tri{\xi^3}{\xi^{15}}{-1}{\xi^{-1}}{\xi}{\xi^{-3}}, \xi\in\G_{54}',\\
	&\tri{\xi^3}{\xi^{-9}}{-1}{\xi^{-1}}{\xi}{\xi^{-3}}, \xi\in\G_{60}',  & \tri{\xi^3}{\xi^{-9}}{-1}{\xi^{-1}}{\xi}{\xi^{-3}}, \xi\in\G_{42}',\\
	& \tri{\xi^3}{\xi^6}{-1}{\xi^{-1}}{\xi}{\xi^{-3}}, \xi\in\G_{42}'.
	\end{align*} 
	\one\,  for $\alpha_1$ and $\alpha_2+\alpha_3$ produces the following lines:
	\begin{align*}
	&   \Dchaintwo{\xi^{3}}{-\xi^{6}}{-1}, \xi\in\G'_{36}  & \Dchaintwo{\xi^{3}}{\xi^{12}}{-1},  \xi\in\G_{54}', && \Dchaintwo{\xi^{3}}{-\xi^{18}}{-1},  \xi\in\G_{60}'\\
	&   \Dchaintwo{\xi^{3}}{-\xi^{9}}{-1}, \xi\in\G'_{42}  & \Dchaintwo{\xi^{3}}{\xi^{3}}{-1},  \xi\in\G_{42}'.
	\end{align*}
	None of these are on the list.
	
	\item In this final case, $\Ft_1=\zeta=q^3$ and $\Ft_2=-1$, which gives $q\in\G_9'\subset \Gf$. 
\end{enumerate}
Hence the lemma follows.
\epf

\begin{lemma}\label{lem:triang-1param-2finite}
	Let $\bq$ be a triangle such that $\bq_{|\I_2}$ is parametric and both $\bq_{|\I_{2,3}}$,  $\bq_{|\{1,3\}}$ are finite. Then $\bq_{|\I_2}=\Gt{q}$ for some $q\in\Gf$ or else $\gdim \B_{\bq}=\infty$.
\end{lemma}
\pf
Let $\Ft=\Dchaintwo{\Ft_1}{\Ft_2}{\Ft_3}$, $\Ft'=\Dchaintwo{\Ft_1'}{\Ft_2'}{\Ft_3'}$, be the diagrams associated to $\bq_{|\I_{2,3}}$ and  $\bq_{|\{1,3\}}$. Then we have a triangle of the form
\[
\tri{\Ft_1'=\Gt{q}_1}{\Gt{q}_2}{\Gt{q}_3=\Ft_1}{\Ft_2}{\Ft_3=\Ft_3'}{\Ft_2'}
\] 
We have that either $q\in\{\Ft_1,\Ft_1'\}$, so $q\in\Gf$, or else $\Gt{q}$ is of type $\sAuno{q}$. Now, if we apply a reflection on vertex 2, we get a triangle:
\[
\tri{q}{q^{-1}}{-1}{\Ft_2^{-1}}{-\Ft_3\Ft_2}{q\Ft_2\Ft_2'}
\]
If $\Dchaintwo{q}{q\Ft_2\Ft_2'}{-\Ft_3\Ft_2}$ is finite, then $q\in\Gf$. If it is parametric, then we are in the setting of the previous lemma, with two parametric diagrams and one finite; hence the lemma follows.
\epf

\begin{lemma}\label{lem:triang-3param}
	Let $\bq$ be a triangle such that all $\bq_{|\I_2}$, $\bq_{|\I_{2,3}}$ and  $\bq_{|\{1,3\}}$ are parametric. Then either $\bq$ belong to the list, or $\bq_{|\I_2}=\Gt{q}$, $\bq_{|\I_{2,3}}=\Gtp{r}$ and $\bq_{|\{1,3\}}=\Gtpp{s}$  for some $q,r,s\in\Gf$ or else $\gdim \B_{\bq}=\infty$.
\end{lemma}

\pf
To start with we can assume that not all vertices are Cartan by Theorem \ref{thm:cartan}. We shall analyze the possible triangles according with the number of vertices labelled with $-1$, from 3 to 0.

\item[$\diamondsuit$] \emph{Three vertices equal to -1.}
Our first case is a triangle
$\tri{-1}{s}{-1}{r}{-1}{q}$.  Note that $q,r,s$ cannot be all simultaneously $-1$, as this would be an example of affine Cartan type. 
We fix $q$ to be the root of 1 of maximal order, hence $q\ne -1$. 
The matrix $\rho_1\bq$ has diagram $\tri{-1}{s^{-1}}{s}{qrs}{q}{q^{-1}}$. We look at the edge $\Dchaintwo{s}{qrs}{q}$. 
If this diagram is not connected (that is, $qrs=1$), then it belongs to rows 9-10-11.
If this is finite, then $q,r,qrs\in\Gf$, so $q,r,s\in \Gf$. Hence we can assume that this edge is parametric: we subdivide into cases, namely when this is of Cartan type, or super type, or other. 

\item[\qquad $\diamond$]Cartan type.

Here $qrs=q^{-\ell}=s^{-1}$ for some $\ell\in\I_3$, so $s=q^\ell$, $r=q^{-2\ell-1}$. The original triangle is thus
$\tri{-1}{q^{\ell}}{-1}{q^{-2\ell-1}}{-1}{q}$, and $\rho_3\bq$ is $\tri{q}{q^{-\ell}}{q^{-2\ell-1}}{q^{2\ell+1}}{-1}{q^{-1}}$. We look at the edge $\lin{q}{q^{-\ell}}{q^{-2\ell-1}}$. If it is finite, then $q\in\Gf$ (hence the same applies to $r$ and $s$). If it is parametric, either $q\in\{\zeta,\zeta^2\}$ or $q^{(2\ell+1)m}q^{-1}=1$ for some $m\in\I_3$. In any case, $q\in\Gf$.

\item[\qquad $\diamond$]Super type.

Assume that $\Dchaintwo{s}{qrs}{q}$ is of super type. Hence $s=-1$ since $\ord s\leq \ord q$ by hypothesis, and $q^{\ell}(-qr)=1$ for some $\ell\in\I_2$, so $r=-q^{-\ell-1}$. Our original triangle becomes $\tri{-1}{-1}{-1}{-q^{-\ell-1}}{-1}{q}$, and $\rho_2\bq$ is $\tri{-1}{-1}{-1}{-q^{\ell+1}}{-q^{-\ell-1}}{q^{-\ell}}$. The tripartite analysis on edge $\lin{-1}{q^{-\ell}}{-q^{-\ell-1}}$ gives that either it is finite, in which case $q^{-\ell},q^{-\ell-1}\in\Gf$ (so $q\in \Gf$), or it is parametric, so there is $m\in\I_2$ such that $(-q^{\ell+1})^mq^{\ell}=1$ and again $q\in\Gf$, or else $\gdim\B_{\bq}=\infty$.

\item[\qquad $\diamond$]Other type.

In this last case edge $\Dchaintwo{q}{qrs}{s}$ is one of the following: $\lin{\zeta}{-1}{-1}$ or $\lin{q}{q^{-1}}{\zeta}$. In the first case $q=\zeta$ and $r=\zeta^2$. As for the second, we have $s=\zeta$, $r=q^{-2}\zeta^2$, so the original triangle and $\rho_2 \bq$ are 
\begin{align*}
&\tri{-1}{\zeta}{-1}{q^{-2}\zeta^2}{-1}{q}, &
&\tri{\zeta}{\zeta^2}{-1}{q^2\zeta}{q^{-2}\zeta^2}{q^{-1}}.
\end{align*}
We perform a tripartite analysis on the edge $\lin{\zeta}{q^{-1}}{q^{-2}\zeta^2}$. If this is finite, then $q\in\Gf$. If it is parametric, then either $q^{-2}\zeta^2=-1$, or $\zeta q^{-1}=1$ or either $q^{-2}\zeta^2q^{-1}=1$, so $q\in\Gf$. 

\item[$\diamondsuit$] \emph{Two vertices labelled with $-1$}. 
Let $t:=q_{22}$ be the vertex such that $t\ne -1$.
Hence either $t=\zeta$, $\widetilde{q}_{2i}=-1$ for some $i\ne 2$, or $2$ is a Cartan vertex.

For the first case we may assume $i=3$. Moreover $s\in\{-1,\zeta,\zeta^2\}$ since the subdiagram with vertices $1$ and $2$ is parametric by hypothesis, so $s\in\Gf$. The diagrams of $\bq$ and $\rho_1\bq$ are, respectively,
\begin{align*}
&\tri{-1}{s}{\zeta}{-1}{-1}{q}, & &\tri{-1}{s^{-1}}{-s\zeta}{-qs}{q}{q^{-1}}. 
\end{align*}
If $q=-s^{-1}$, then $q\in\Gf$. Otherwise $-qs\ne 1$, and 
$\rho_1\bq$ is a triangle containing the subdiagram $\lin{-s\zeta}{-qs}{q}$. If this subdiagram does not belong to the list, then $\gdim\B_{\bq}=\infty$. Otherwise it is either parametric, which says that $q\in\{-1,\zeta,\zeta^2\}$ or $q^c(-qs)=1$ for some $c\in\I_3$, or else finite; in any case $q\in\Gf$.

Now assume that $t$ is a Cartan vertex: the diagram of $\bq$ is $\tri{-1}{t^{-a}}{t}{t^{-b}}{-1}{q}$ for some $a,b\in\I_3$. \one\,  for $\alpha_1+\alpha_3$, $\alpha_2$ gives the diagram $D:=\lin{q}{t^{-a-b}}{t}$: if $\mathtt{d}$ is not in the list, then $\gdim\B_{\bq}=\infty$. If this is finite, then $q,t\in\Gf$. If $d$ is not connected, i.e. $t^{a+b}=1$, then $t\in\Gf$ and $\rho_1\bq$ contains the diagram $\lin{q}{q}{-t^{1-a}}$: a tripartite analysis ends this case. Otherwise, $d$ is parametric. We take care of different choices of $a,b$:
\begin{itemize}[leftmargin=*]
	\item $a=3$ or $b=3$. We may fix $a=3$: this forces $\bq_{|\I_2}$ to be of Cartan type $G_2$ with $t\in\G_6'$. We get a contradiction since $t\ne -1,\zeta,\zeta$ and $t^m(-t^{-b})\ne 1$ for $m=1,2,3$. 
	\item $a=b=2$. As $t\ne \pm 1$, either $t^3=1$, $qt^{-4}=1$, or else $t^{m-4}=1$ for some $m\in\I_3$: clearly $m=2,3$ are not allowed, so $t\in\G_3'$ and $q\in\{-1,t^{\pm1}\}$.
	\item $\{a,b\}=\{1,2\}$. As $\ord t\ge 4$ (since we assume that both $\bq$ and $\mathtt{d}$ are connected), vertex $2$ of $\mathtt{d}$ is of Cartan type: $t^{m-3}=1$ for some $m\in\I_3$. This forces $m=3$, which implies as well $q=t^3$. Hence $\bq$ belongs to row 7 of the list.
	\item $a=b=1$. As $t=-1$ by hypothesis, either $t^3=1$, $qt^{-2}=1$ or else $t^{m-2}=1$ for some $m=1,2,3$: only $m=2$ is allowed, which forces $q=t^2$ or $q=-1$. For $q=t^2$ we have row 6 of the list, while for $t=-1$, $\rho_3 \bq$ is a triangle with $-1$ in the three vertices and was already considered.
\end{itemize}

\item[$\diamondsuit$] \emph{One vertex labelled with $-1$}.
We fix $q_{11}=-1$ and start by distinguishing three possible cases, namely
\begin{itemize}[leftmargin=*]
	\item[(a)] $q_{jj}=\zeta\in\G_3'$, $\qti_{1j}=-1$, $j=2,3$.
	\item[(b)] $q_{jj}=\zeta\in\G_3'$ for a single $j$, say $j=2$, $\qti_{1j}=-1$ and $q=q_{33}\neq -1,\zeta$, so $q_{33}^m\qti_{13}=1$, $m\leq 3$.
	\item[(c)] $q_{22}^a\qti_{12}=1=q_{33}^b\qti_{13}$, $a,b\leq 3$.
\end{itemize}

We analyze each case. For (a), we get the triangle $\tri{-1}{-1}{\zeta}{s}{\zeta}{-1}$, for which \one\, for $\alpha_1+\alpha_2,\alpha_3$ gives $s=\zeta^2$.

As for (b), we have $\tri{-1}{-1}{\zeta}{r}{q}{q^{-m}}$, and we focus on the line $\lin{\zeta}{r}{q}$. If it is of Cartan type, then it follows that $q,r\in\Gf$. If not, then we necessarily have $r=q^{-1}$. Back to the triangle, we apply \one\, for $\alpha_1+\alpha_2,\alpha_3$ and obtain $\lin{\zeta}{q^{-1-m}}{q}$. This diagram is in the list if and only if it is disconnected ($q^{1+m}=1$), is finite (so $q\in\Gf$) or it is parametric, which gives $q^dq^{-1-m}=1$, for some $d\leq 3$, so $q\in\Gf$.

For case (c), we consider two possibilities, namely
\begin{itemize}[leftmargin=*]
	\item[(c1)] $r=\zeta\in\G_3'$ and $qs=1$ (or, equivalently, $q=\zeta$, $rs=1$).
	\item[(c2)] $\lin{r}{s}{q}$ is of Cartan type, so we may assume $r=q^m$, $s=q^{-m}$, $m\leq 3$.
\end{itemize}

For case (c1), \one\, with $\alpha_1+\alpha_3, \alpha_2$ generates the line $\lin{-\zeta^{x}}{q^{y}}{q}$, $x=-1-b$, $y=-1-a$. If this diagram is either disconnected or finite, then the case follows. If it is parametric, then we look at the possible values of $a$. If $a=3$, then $q\in\G_6'$ ($\lin{-1}{q^{-3}}{q}$, $q^3=-1$) is of type $\Gdost{q}$. If $a=2$, then the line $\lin{-\zeta^{x}}{q^{y}}{q}=\lin{-\zeta^{x}}{q^{-3}}{q}$ is again $\Gdost{q}$ and thus $q^3=-\zeta^{-1-b}$ gives $\ord q\leq 30$. If $a=1$, the reflection $\rho_1$ generates a triangle with at least two vertices with parameter -1 or a line, cases already analyzed.

In case (c2), we consider the triangle $\tri{-1}{q^{-a}}{q}{q^{-m}}{q^m}{q^{-bm}}$. We analyze the possible values of $(a,b)$. If $a=b=1$, then once again  $\rho_1$ leads to a triangle with at least two vertices with parameter -1 or to a line. As above, if $a=3$, then $\ord q=6$ (type $\Gdost{q}$).  Similarly, if $b=3$, then $\ord q^m=6$ and then $q\in\Gf$. Finally, assume $a=b=2$, that is $\tri{-1}{q^{-2}}{q}{q^{-m}}{q^m}{q^{-2m}}$. In this last case, \one\, for $\alpha_1,\alpha_2+\alpha_3$ gives $\lin{-1\quad}{q^{-2m-2}}{q}$. If this is either disconnected or finite, the case follows. If it is parametric, then either $q^{3}=1$ or there is $d\in\I_3$ such that $q^dq^{-2m-2}=1$. In any case, $q\in\Gf$.

\item[$\diamondsuit$] \emph{No vertices with label $-1$}.

We fix, without loss of generality, $q_{11}=\zeta$, $q=q_{22}\notin \G_3'$ with $\qti_{12}q=1$, so our triangle is of the form $\tri{\zeta}{q^{-1}}{q}{r}{s}{t}$.
Again, we study some subcases, according to (a) $st=1$, (b) $s^2t=1$,  or (c) $s^3t=1$.

In cases (b) or (c) we also have  $t\zeta=1$, so $s\in\G_6'$ as above. Hence, $\lin{s}{r}{q}$ is of Cartan type with $q,r\in\Gf$.

For case (a), we look at edge $\lin{s}{r}{q}$. We either have (a1) $s^3=1$ and $qr=1$ or (a2) $s\notin\G_3'$ so $s=q^m=r^{-1}$ for some $m\in\I_3$. In (a1), \one\, for $\alpha_1+\alpha_3,\alpha_2$ gives $\lin{\zeta}{q^{-2}}{q}$, so $\ord q\leq 6$.¸ In (a2), \one\, for $\alpha_1+\alpha_3,\alpha_2$ gives $\lin{\zeta}{q^{-1-m}}{q}$ and therefore $q\in\Gf$.

This ends all possible triangles as in the lemma.
\epf

\subsection{\gap\ lemmas}\label{subsec:gap}

In this section we use \gap\, to work with all possible lines and triangles made up of gluing together diagrams of rank 2, either finite or parametric, these last ones evaluated at roots in $\Gf$. We consider the criteria developed in \S \ref{sec:crit} to rule out each such diagram, we also discard those already belonging to the list. This process is carried out with the help of computer software \gap, in two lemmas, one dealing with lines and another one dealing with triangles. On each proof, we cite the code files (\tt{lines/triangles.g}) used and the log files obtained (\tt{lines/triangles.log}). In turn, these code files rely on definitions stated in \tt{basic.g}. All of these files are stored in the authors' webpages.

\begin{lemma}\label{lem:gaplines}
	Assume that $d$ is a line obtained by pasting two diagrams, either finite or parametric evaluated at a root in $\Gf$.
	Then $\gdim\B_\bq=\infty$ or $d$ belongs to the list.
\end{lemma}
\pf
We check with \gap\, all possible combinations. See \tt{lines.(g|log)}.  After all criteria, we have the following 7 diagrams left:
\begin{align*}
\liin{         -1}{               -1}{          -\zeta^2}{            -\zeta}{             \zeta };&&
\liin{         -1}{            -\zeta}{             \zeta}{            -\zeta}{          -\zeta^2 };&&
\liin{         -1}{            -\zeta}{             \zeta}{          -\zeta^2}{            -\zeta };\\
\liin{     \zeta^2}{             \zeta}{               -1}{               -1}{             \zeta };&&
\liin{       \zeta}{           \zeta^2}{               -1}{            -\zeta}{             \zeta };&&
\liin{       \zeta}{            -\zeta}{          -\zeta^2}{           \zeta^2}{             \zeta }\\
&& \liin{      -\xi^4}{          \xi^2}{        -\xi^{-1}}{          -\xi}{             \xi^4 };&&  
\end{align*}
for $\zeta\in\G_3',  \xi\in\G_{12}'$. Now, the last one, involving $\xi\in\G_{12}'$, is standard with affine matrix $A^{(2)}_4$; the same is true for the sixth diagram. The first is standard with indefinite Cartan matrix. The second and third are of affine type $D_3^{(2)}$. We recall the diagrams $\bq^{(n)}$, $n\in\I_5$, from Example \ref{ex:basic-datum-6diagrams}: the fourth diagram above is $\bq^{(5)}$ up to $\zeta\leftrightarrow\zeta^2$ and the fifth is $\bq^{(4)}$ up to permutation; hence Lemma \ref{lem:grupoide6puntos} applies. Thus the Nichols algebras of these seven diagrams have infinite $\gdim$.
\epf

\begin{lemma}\label{lem:gaptriangles}
	Assume that $d$ is a triangle obtained by pasting three diagrams, either finite or parametric evaluated at a root in $\Gf$.
	Then $\gdim\B_\bq=\infty$ or $d$ belongs to the list.
\end{lemma}
\pf
Again, we check with \gap\, all possible combinations of these three diagrams, either finite or parametric, evaluated at a root in $\Gf$; we also discard those triangles than can be reflected to a line, as we have already dealt with that case in the previous lemma. See \tt{triangles.(g|log)}, where we leave out the case of three  edges of type $\sAuno{q}$, analyzed below.  

After all criteria, we have the following 5 diagrams left:
\begin{align*}
\tri{         -1}{       -1}{         -1}{  -\zeta}{         -1}{  -\zeta};&&
\tri{         -1}{       -1}{         -1}{  -\zeta}{    -\zeta^2}{  -\zeta };&&
\tri{    -\zeta^2}{  -\zeta}{         -1}{  -\zeta}{    -\zeta^2}{  -\zeta };\\
\tri{  -\xi_4}{   \xi_4}{     -1}{   \xi_4}{  -\xi_4}{   \xi_4 };&&
\tri{   \xi_4}{  -\xi_4}{     -1}{  -\xi_4}{   \xi_4}{  -\xi_4 }
\end{align*}
On the one hand, notice that triangles 4 and 5 are actually equivalent, by changing the root $\xi\leftrightarrow-\xi$ in $\G_4'$. Also,  $\tau_{23}\rho_2$ of triangle one gives triangle three, with $\zeta\leftrightarrow\zeta^2$. Similarly,  $\rho_1$ of the second triangle gives the first triangle, again with $\zeta\leftrightarrow\zeta^2$.  Hence we need to deal with matrices $\bq$ and $\bq'$ with diagrams 1 and 4, respectively. The corresponding basic data are
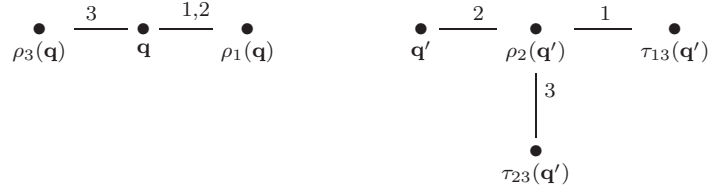
\begin{align*}
&\xymatrix@C-4pt{\underset{\rho_3(\bq)}{\bullet} \ar @{-}[r]^{3} &
	\underset{\bq}{\bullet} \ar @{-}[r]^{1,2} & \underset{\rho_1(\bq)}{\bullet} }
&
&\xymatrix@C-4pt{\underset{\bq'}{\bullet} \ar @{-}[r]^{2} &
	\underset{\rho_2(\bq')}{\bullet} \ar @{-}[r]^{1} \ar @{-}[d]^{3} & \underset{\tau_{13}(\bq')}{\bullet} 
	\\
	& \underset{\tau_{23}(\bq')}{\bullet} & }
\end{align*}
Both are standard with Cartan matrix $A_2^{(1)}$, and Proposition \ref{prop:standard-affine} applies.

Let us now focus on merging $\sAuno{q},\sAuno{r},\sAuno{s}$. That is
\[
\tri{-1}{q}{-1}{r}{-1}{s}, \qquad q,r,s\in\k\setminus\{\pm1\}.\]
Notice that when $qrs=1$, then it is of finite type (and it can be reflected to a line, which is true also when  $qr=1, qs=1$ or $rs=1$). In any case, a reflection on any vertex returns a triangle of a different type or forces $q,r$ or $s$ to be -1, a contradiction.
\epf


\begin{thebibliography}{AAG1}
	
	\bibitem[A]{A-leyva} N. Andruskiewitsch. \emph{An Introduction to Nichols Algebras}. In Quantization, Geometry and Noncommutative Structures in Mathematics and Physics. 
	A. Cardona, P. Morales, H. Ocampo, S. Paycha, A. Reyes, eds., pp. 135--195, Springer (2017).
	
	
	\bibitem[AA1]{AA}
	N. Andruskiewitsch, I. Angiono. \emph{On Nichols algebras with generic braiding},
	in  \emph{Modules and
		Comodules}, T. Brzezinski; J. L. G\'omez Pardo;
	I. Shestakov; P. F. Smith (Eds.). Trends in Mathematics (2008), 47--64.
	
	
	\bibitem[AA2]{AA-diag-survey} \bysame \emph{On Finite dimensional Nichols algebras of diagonal type}. 
	Bull. Math. Sci. \textbf{7} (2017), 353--573.
	
	
	\bibitem[AAH1]{AAH-conj} 
	N. Andruskiewitsch, I. Angiono, I. Heckenberger, \emph{On finite GK-dimensional Nichols algebras over abelian groups}.  Mem.~AMS, to appear.
	
	
	\bibitem[AAH2]{AAH-rank2} 
	\bysame \emph{On finite GK-dimensional Nichols algebras of diagonal type}.  
	Contemp. Math. \textbf{728} (2019), 1--23.
	
	\bibitem[AAY]{AAY}
	N. Andruskiewitsch, I. Angiono, M. Yakimov, \emph{Poisson orders on large quantum groups}, \texttt{arXiv:2008.11025}.
	
	\bibitem[ASa]{ASa} N. Andruskiewitsch, G. Sanmarco.
	\emph{Finite GK-dimensional pre-Nichols algebras of quantum linear spaces and of Cartan type}. Trans. AMS Ser. B \textbf{8}, 296--329 (2021).
	
	
	\bibitem[ASc]{AS-adv}  N. Andruskiewitsch, H.-J. Schneider, 
	\emph{Finite quantum groups and Cartan matrices}. Adv. Math. \textbf{154} (2000), no. 1, 1--45. 
	
	\bibitem[ACS]{ACS} I. Angiono, E. Campagnolo, G. Sanmarco.
	\emph{Finite GK-dimensional pre-Nichols algebras of super and standard type}, \texttt{arXiv:2009.04863}.
	
	\bibitem[C+]{C} L. Carbone, S. Chung, L. Cobbs, R. McRae, D. Nandi, Y. Naqvi, D. Penta, \emph{Classification of hyperbolic Dynkin diagrams, root lengths and Weyl group orbits}. 
	J. Phys. A \textbf{43} (2010), no. 15, 155209, 30 pp.
	
	\bibitem[H1]{H-inv} I. Heckenberger, 
	{\em The Weyl groupoid of a Nichols algebra of diagonal type}, Inventiones Math. \textbf{164} (2006), 175--188.
	
	\bibitem[H2]{H-rk2} \bysame \emph{Rank 2 Nichols algebras with finite arithmetic root system}. Algebr. Represent. Theory \textbf{11} (2008), no. 2, 115--132.
	
	\bibitem[H3]{H-rk3} \bysame \emph{Classification of arithmetic root systems of rank 3}. Proceedings of the XVIth Latin American Algebra Colloquium (Spanish), 227--252, Bibl. Rev. Mat. Iberoamericana, Rev. Mat. Iberoamericana, Madrid, 2007.
	
	\bibitem[H4]{H-full} \bysame \emph{Classification of arithmetic root systems},
	Adv. Math. \textbf{220},  59--124 (2009).
	
	\bibitem[HS]{HS} I. Heckenberger, H.-J. Schneider, 
	\emph{Root systems and Weyl groupoids for Nichols algebras}, Proc. London Math. Soc. \textbf{101} (2010), 623--654.
	
	\bibitem[K]{K} V. Kac, \emph{Infinite-dimensional Lie algebras}, 3rd ed., Cambridge University Press, Cambridge, 1990.
	
	\bibitem[Kh1]{K1} V. Kharchenko, \emph{A quantum analog of the Poincare-Birkhoff-Witt theorem}, Algebra	and Logic \textbf{38} (1999), 259--276.
	
	\bibitem[Kh2]{K2} \bysame \emph{PBW-bases of coideal subalgebras and a freeness theorem}
	Trans.~AMS \textbf{360} (10) (2008), 5121--5143.
	
	\bibitem[KL]{KL} G. Krause, T. Lenagan, \emph{Growth of algebras and Gelfand-Kirillov dimension}. Revised edition.
	Graduate Studies in Mathematics, \textbf{22}. American Mathematical Society, Providence, RI, 2000. x+212 pp
	
	\bibitem[R]{R quantum groups} M. Rosso. \emph{Quantum groups and quantum shuffles}. Invent. Math. \textbf{133} (1998), 399--416.
	
\end{thebibliography}
\end{document}